\documentclass[leqno]{amsart}
\usepackage{amsmath,amssymb}

\usepackage[all]{xy}

\newcommand{\Jac}{\mathop{\mbox{\rm Jac}}\nolimits}

\begin{document}
\begin{center}
\title{Hyperelliptic $d$-osculating covers and rational surfaces}
\author{Armando Treibich}
\maketitle
\end{center}
  
\noindent \medskip\section{ Introduction}

\noindent\textbf{1.1.}
Let $\mathbb{P}^1$ and $(X,q)$ denote, respectively, the projective line and a fixed  
elliptic curve marked at its origin, both defined over an algebraically closed field  $\mathbb{K}$ of arbitrary characteristic $\emph{\textbf{p}} \neq2$.
We will study all finite separable marked morphisms $\pi :(\Gamma ,p)\rightarrow (X,q)$, called hereafter \textit{hyperelliptic covers}, such that $\Gamma$ is a degree-$2$ cover of  $ \mathbb{P}^1$, ramified at the smooth point $p \in \Gamma$. Canonically associated to $\pi$ there is the Abel (rational) embedding of $\Gamma$ into its \emph{generalized Jacobian}, $A_p: \Gamma \to Jac\,\Gamma$, and $\{0\} \subsetneq  V^1_{\Gamma,p}\ldots \subsetneq V ^g_{\Gamma,p}$, the flag of hyperosculating planes to $A_p(\Gamma)$ at $A_p(p)\in Jac\,\Gamma$ (cf. \textbf{2.1. \& 2.2.}).
On the other hand, we also have the homomorphism $\iota_\pi: X \to \Jac\,\Gamma$, obtained by dualizing $\pi$. There is a smallest positive integer $d$ such that the tangent line to $\iota_\pi( X)$ is contained in the $d$-dimensional osculating plane $V^d_{\Gamma,p}$. We call it \emph{the osculating order} of $\pi$, and $\pi$ a \emph{hyperelliptic $d$-osculating cover} (\textbf{2.4.}(2)). If $\pi$ factors through another \emph{hyperelliptic cover}, the arithmetic genus increases, while the \emph{osculating order} can not decrease (\textbf{2.8.}).
\indent Studying, characterizing and constructing those with given \textit{osculating order} $d$ but maximal possible arithmetic genus, so-called \emph{minimal-hyperelliptic d-osculating covers}, will be one of the main issues of this article. The other one, to which the first issue reduces, is the construction of all rational curves in a particular anticanonical rational surface associated to $X$ (i.e.: a rational surface with an effective anticanonical divisor). Both problems are interesting on their own and in any characteristic $\emph{\textbf{p}} \neq2$. They were first considered, however, over the complex numbers and through their link with solutions of the \textit{Korteweg-deVries} hierarchy, doubly periodic with respect to the $d$-th \textit{KdV} flow  (cf. \cite{A-McK-M}, \cite{D-N}, \cite{I-M}, \cite{K}, \cite{T-V} for $d=1$ and \cite {S}, \cite{A-K-V}, \cite{F}, \cite{F-T} for $d = 2$). 
We sketch hereafter the structure and main results of our article.

  \begin {enumerate}
  
 \item We start defining in section \textbf{2.} the Abel rational embedding $A_p:\Gamma \to Jac\,\Gamma$,  and construct the flag  $\{0\} \subsetneq  V^1_{\Gamma,p}\ldots \subsetneq V^g_{\Gamma,p} = H^1(\Gamma,O_\Gamma)$, of \emph{hyperosculating planes} at the image of any smooth point $p \in \Gamma$. We then define the homomorphism $\iota_\pi : X \to Jac\,\Gamma$, canonically associated to the \emph{hyperelliptic cover} $\pi$, and its \emph{osculating order} (\textbf{2.4.}(2)). Regardless of the \emph{osculating order}, we prove that any degree-$n$ \emph{hyperelliptic cover} has odd ramification index at the marked point, say $\rho$, and factors through a unique one of maximal arithmetic genus $2n\,$-$\,\frac{\rho+1}{2}$ (\textbf{2.6.}). We finish characterizing the \emph{osculating order} by the existence of a particular projection $\kappa:\Gamma \rightarrow \mathbb{P}^1$ (\textbf{2.6.}). \medskip 
 
\item The $d$-osculating criterion \textbf{2.6.} paves the way to the algebraic surface approach
  developed in the remaining sections. The main characters are played by (two morphisms between)
three projective surfaces, canonically associated to the elliptic curve $X$:
   
   \medskip $\centerdot \,\,e:S^\bot \to S$ : the blowing-up of a particular ruled surface $\pi_S : S \rightarrow X$, at the $8$ fixed points of its involution; 
   
   \medskip  $\centerdot \,\,\varphi: S^\bot \to \widetilde S$ : \,\,a projection onto an anticanonical rational surface. \medskip
  
\item Once $S$, $S^\bot$ and $\widetilde S$ are constructed (\textbf{3.2.}, \textbf{3.4.}), we prove that any \textit{hyperelliptic 
  $d$-osculating cover} $\pi: (\Gamma, p) \to (X,q)$ factors canonically through a curve $\Gamma^\bot \subset S^\bot$, and projects, via $\varphi : S^\bot \to \widetilde S$, onto a rational irreducible curve $\widetilde \Gamma \subset \widetilde S$  (\textbf{3.8.}). We also prove that any \textit{hyperelliptic $d$-osculating cover} dominates a unique one of same \textit{osculating order} $d$, but maximal arithmetic genus, so-called \textit{minimal-hyperelliptic} (\textbf{3.9.}). Conversely, given $\widetilde \Gamma \subset \widetilde S$, we study when and how one can recover all \textit{minimal-hyperelliptic $d$-osculating covers} having same canonical projection $\widetilde \Gamma $ (\textbf{3.11.}) .\\
      
\item  Section \textbf{4.} is mainly devoted to studying the linear equivalence class of the curve $ \Gamma^\bot \subset S^\bot$, canonically associated to any \textit{hyperelliptic $d$-osculating cover} $\pi$, and associated invariants  (\textbf{4.3.} \& \textbf{4.4.}). We end up with a numerical characterization of \textit{minimal-hyperelliptic $d$-osculating covers} (\textbf{4.6.}).\medskip

\item At last, we dress the list of all $(\,$-$\,1)$ and $(\,$-$\,2)$-irreducible curves of $\widetilde S$ (\textbf{5.7.}), needed to study its \textit{nef} cone, and give, for any $n,d \in \mathbb{N}^*$, two different constructions of  $(d\,$-$\,1)$-dimensional families of smooth, degree-$n$, \emph{minimal-hyperelliptic $d$-osculating covers}: one based on Brian Harbourne's results on anticanonical rational surfaces (\cite{Har}), the other one based on \cite{T.2} and leading, ultimately, to explicit equations for the corresponding covers. \\
      
  \end{enumerate}

\noindent  \section{ Jacobians of curves and hyperelliptic $d$-osculating covers}

 \noindent\textbf{2.1.}
	 	 Let $\mathbb K$ be an algebraically closed field of characteristic $\emph{\textbf{p}}\neq 2$, $\mathbb
P^1$ the projective line  over $\mathbb K$ and $(X,q)$ a fixed elliptic curve, also defined over $\mathbb K$. The latter will be equipped with its canonical symmetry $[\,$-$\,1]: (X,q) \to (X,q)$, fixing $\omega_o:=q$, as well as the other three half-periods $\{\omega_j, \, j=1,2,3\}$. We will also choose once for all, an odd local parameter of $X$ centered at $q$, say $z$, such that $z \circ [\,$-$\,1]=\,$-$\,z$.  

By a curve we will
mean hereafter a complete integral curve over $\mathbb K$, say $\Gamma$, of positive
 arithmetic genus $ g>0$. The moduli space of degree-$0$ invertible sheaves over $\Gamma$, denoted by $Jac\,\Gamma$ and called the \textit{generalized Jacobian} of $\Gamma$, is a $g$-dimensional connected commutative algebraic group, canonically identified to  $H^1(\Gamma,O_\Gamma ^{*})$, with tangent space at its origin equal to $H^1(\Gamma,O_\Gamma)$. 	
 Recall also the \textit{Abel} (rational) embedding $A_p : \Gamma \rightarrow Jac\,\Gamma$, sending any smooth point $p'\in \Gamma $ to the isomorphism class of $O_\Gamma(p'\,$-$\,p)$. For
any marked curve $(\Gamma , p)$ as above, and any positive integer
$j$, let us consider the exact sequence of $O_\Gamma $-modules $0\to
O _\Gamma \to O_\Gamma (jp)\to O_{jp} (jp)\to 0$, as well as the
corresponding long exact cohomology sequence :\\

$0 \to H^0(\Gamma,O_\Gamma ) \to H^0\big(\Gamma,O_\Gamma (jp)\big)\to H^0\big(\Gamma,O_{jp}(jp)\big)\stackrel{\delta }\to H^1(\Gamma,O_\Gamma )\to \ldots ,$\\

where $\delta : H^0\big(\Gamma,O_{jp} (jp)\big) \to H^1(\Gamma,O_\Gamma)$ is the cobord morphism. 
According to the Weierstrass gap Theorem, for any
$d \in\{1,\ldots ,g \}$, there
exists $0<j < 2g$ such that $\delta \Big(H^0\big(\Gamma,O_{jp}(jp)\big)\Big)$ is a
$d$-dimensional subpace, denoted hereafter by $V^d_{\Gamma,p}$. 

For a
generic point $p$ of $\Gamma $ we have $V^d_{\Gamma,p}= \delta
\Big(H^0\big(\Gamma,O_{dp}(dp)\big)\Big)(i.e.:j=d)$,  while for any $p\in \Gamma$, the tangent to $A_p(\Gamma )$ at $0$ is equal to $V^1_{\Gamma,p}=\delta \Big(H^0\big(\Gamma,O_p(p)\big)\Big)$.\\

\textbf{Definition 2.2.}\\
\indent (1) \emph{ The filtration
$\{0\} \subsetneq  V^1_{\Gamma,p}\ldots \subsetneq V^g_{\Gamma,p} = H^1(\Gamma,O_
\Gamma)$ will be called  the flag of hyperosculating spaces to $A_p(\Gamma )$ at }$
0 $.\\

\indent (2)  \emph{The curve $\Gamma$ will be called a hyperelliptic curve, and $p\in \Gamma$ a Weierstrass point, if there exists a degree-$2$ projection onto $\mathbb{P}^1$, ramified at $p$. Or equivalently, if  there exists an involution, denoted in the sequel by $\tau_\Gamma: \Gamma \to \Gamma$ and called the hyperelliptic involution, fixing $p$ and such that the quotient curve $\Gamma/\tau_\Gamma$ is isomorphic to $\mathbb{P}^1$.}\\

 \textbf{Proposition 2.3. }(\cite{T.1}\S1.6.)
 
\textit{Let $(\Gamma, p, \lambda) $ be a hyperelliptic curve of arithmetic genus $g$, equipped with a local parameter $\lambda$, centered at a smooth 
Weierstrass point $p \in \Gamma$. For any odd integer $1\leq j:= 2d\,$-$\,1\leq g$, consider the exact sequence of $O_\Gamma$-modules}: 
\begin{displaymath}
 0\to O_\Gamma \to O_\Gamma (jp)\to O_{jp}(jp) \to 0  \quad , 
\end{displaymath}
\textit{as well as its long exact cohomology sequence}
 
\begin{displaymath}
 0\to H^0(\Gamma,O_\Gamma )\to H^0\big(\Gamma,O_\Gamma (jp)\big)\to H^0\big(\Gamma,O_{jp}(jp)\big) \stackrel {\delta }{\to }  H^1(\Gamma,O_\Gamma) \to \dots,
\end{displaymath}

$\delta$  \textit{   being the cobord morphism.}\\

\noindent \textit{For any, $m\geq 1$, we also let $ \, [\lambda^{\textrm{-}\,m}]$ denote the class of $\lambda^{\textrm{-}\,m}$ in $H^0\big(\Gamma,O_{mp}(mp)\big)$.
Then $V^d_{\Gamma,p}$ is generated by $\Big\{ \delta \big([\lambda^{2l\textrm{-}1}]\big),\, l = 1,..,d \Big\}$. In other words, the $d$-th osculating subspace to $A_p(\Gamma )$ at $ 0$ is equal to $\delta\Big(H^0\big(\Gamma,O_{jp}(jp)\big)\Big)$, for} $j=2d\,$-$\,1$.\\

\textbf{Definition 2.4.} 

\indent (1) \textit{A finite separable marked morphism $\pi:(\Gamma,p)\to(X,q)$, such that $\Gamma$ is a hyperelliptic curve and $p \in \Gamma$ a smooth Weierstrass point, will be called a hyperelliptic cover. We will say that $\pi$ dominates another hyperelliptic cover $\overline \pi:(\overline \Gamma,\overline p)\to(X,q)$, if there exists a degree-$1$ morphism $j:(\Gamma, p) \to (\overline \Gamma, \overline p)$, such that  $\pi = \overline \pi \circ j$.}\\

\indent (2) \textit{Let
$\iota_\pi:X\to \ Jac\,\Gamma $ denote the group
homomorphism $q'\mapsto A_p\big(\pi^*(q'\textrm{-}\,\,q)\big)$. There is a minimal integer $d\geq 1$, called henceforth osculating order of $\pi$, such that the tangent to $\iota_\pi(X)$ at $ 0 $ is
contained in $V^d_{\Gamma,p}$. We will then call $\pi$ a hyperelliptic $d$-osculating cover.} \\

 \textbf{Proposition 2.5.}
 
\noindent \textit{Let $\pi :(\Gamma,p)\to(X,q)$ be a degree-$n$ hyperelliptic cover with ramification index $\rho$ at $p$, }$f: (\Gamma, p) \to (\mathbb{P}^1, \infty)$ \textit{the corresponding degree-$2$ projection, ramified at $p$, and let $\Gamma_{f,\pi}$ denote the image curve $(f,\pi)(\Gamma)\subset \mathbb{P}^1 \times X$. Then (see diagram below)}, 
 
 \begin{enumerate}
 \item  \textit {the hyperelliptic involution $\tau_\Gamma$ satisfies$\quad[\,$-$\,1] \circ \pi=\pi \circ \tau_\Gamma$ and $\rho$ is odd};
 
 \item \textit{$\Gamma_{f,\pi}$ has arithmetic genus $2n\,$-$\,1$ and is unibranch at $(\infty,q)$};
 
 \item \textit{let $(\overline \Gamma, \overline p)$ denote the partial desingularization of $\Gamma_{f,\pi}$ at $(\infty,q)$, equipped with its canonical projection via $\Gamma_{f,\pi}$, say $\overline \pi: (\overline \Gamma, \overline p) \to (X,q)$,  then:}\\

$\overline \pi$ \textit{is a hyperelliptic cover of arithmetic genus} $2n\,$-$\,\frac{1}{2}(\rho+1)$;\\
 
 \item \textit{$\pi$, as well as any hyperelliptic cover dominated by $\pi$, factors through} $\overline \pi$.
 \end{enumerate}

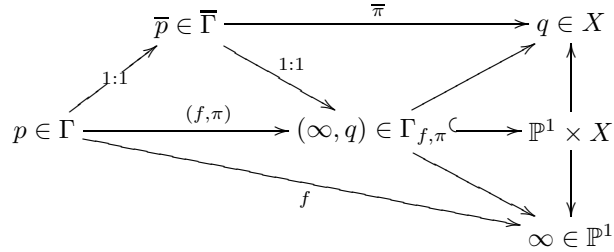
 \begin{displaymath}
\xymatrix{
&\overline p\in \overline \Gamma \ar[rr]^{\overline \pi}\ar[dr]^{1:1}  &  &  q \in X\\
p\in\Gamma \ar[ru]|{1:1}\ar[drrr]_{f}\ar[rr]^{(f,\pi)} & &(\infty,q) \in \Gamma_{f,\pi} \ar[ur] \ar[dr]\ar@{^{(}->}[r]&  \mathbb{P}^1\times X \ar[u] \ar[d]  \\
& && \infty \in \mathbb{P}^1 
     }
\end{displaymath}

 \textbf{Proof.} (1) Let  $Alb_\pi : Jac\,\Gamma \to Jac\,X $ denote the Albanese homomorphism, sending any $L\in Jac\,\Gamma$ to $ Alb_\pi(L):= det(\pi_*L)\otimes det(\pi_*O_\Gamma)^{-1}$, and $\Gamma^0 $ denote the open dense subset of smooth points of 
$\Gamma$. Up to identifying $Jac\,X$ with $(X,q)$, we know that $Alb_\pi \circ \iota_\pi = [n]$, the multiplication by $n$, and $Alb_\pi \circ A_p$ is well defined over $\Gamma^0$ and equal there to $\pi$. Knowing, on the other hand, that $A_p \circ \tau_\Gamma = [-1] \circ A_p$, we deduce that $\pi \circ \tau_\Gamma = Alb_\pi \circ A_p \circ \tau_\Gamma= [-1] \circ Alb_\pi \circ  A_p= [-1] \circ \pi$ (over the open dense subset $\Gamma^0$, hence) over all $\Gamma$ as asserted.

(2) \& (3) The projections $f$ and $\pi$ have degrees $2$ and $n$, implying that $\Gamma_{f,\pi}$ is numerically equivalent to $n. \{\infty \}\times X + 2 . \mathbb{P}^1 \times \{q\}$ and, by means of the adjunction formula, that it has arithmetic genus $2n\,$-$\,1$. We also know that $f$ and $\pi$ have ramification indices $2$ and $\rho$ at $p \in \Gamma$. Hence,  $\Gamma_{f,\pi}$ intersects the fibers $\mathbb{P}^1\times\{q\}$ and $\{\infty\}\times X$ at $(\infty,q)$, with multiplicities $\rho$ and $2$. Adding property \textbf{2.5.}(1) we deduce that its local equation at $(\infty,q)$  can only have even powers of $z$, and must be equal to $z^2 = w^\rho h(w,z^2)$, for some invertible element $h$ (i.e.: $h(0,0)\neq0)$. In particular $\Gamma_{f,\pi}$  is unibranch and has multiplicity $min\{2, \rho\} $ at $(\infty,q)$. Moreover, for its desingularization over $(\infty,q)$, $\frac{\rho -1}{2}$ successive monoidal transformation are necessary, each one of which decreases the arithmetic genus by $1$. Hence $\overline \Gamma$ has arithmetic genus $2n\,$-$\,1\,$-$\,\frac{\rho -1}{2}=2n\,$-$\,\frac{\rho +1}{2}$ as asserted.

(4) Since $\Gamma$ is already smooth at $p$, we immediately see that $(f,\pi)$ factors through $\overline \pi$. Hence, $\pi$ dominates $\overline \pi$ as asserted. Reciprocally, any other \textit{hyperelliptic cover} dominated by $\pi$ must factor through $\big(\Gamma_{f,\pi}, (\infty,q)\big)$, and should lift to its partial desingularization $(\overline \Gamma, \overline p)$. In other words, it should dominate $\overline \pi$.$\blacksquare$\\ 
 
\textbf{ Theorem 2.6.}

 \noindent \emph{The osculating order of an hyperelliptic cover  $\pi :(\Gamma ,p)\rightarrow (X,q)$, is the minimal integer $d\geq1$ for which there exists a morphism $
\kappa :\Gamma \rightarrow \mathbb{P}^1$ satisfying :} \\

(1) \emph{the poles of $\kappa $ lie along $\pi ^{\textrm{-}1}(q)$};

(2) $\kappa $+$\,\pi ^{*}(z^{\,\textrm{-}1})$ \emph{has a pole of order $2d\,$-$1$  at $p$, and no
other pole along $\pi ^{\textrm{-}1}(q)$} (\textbf{2.1.}).\\

\noindent \emph{Furthermore, for such $d$ there exists a unique morphism} $
\kappa :\Gamma \rightarrow \mathbb{P}^1$ \emph{satisfying properties} (1)\&(2)
\emph{above, as well as }(\textbf{2.2.}(2)):

\medskip\ (3) $\tau _\Gamma ^{*}(\kappa )=\,$-$\,\kappa $.\\

\textbf{Proof.} According to \textbf{2.3.}, $ \forall r\geq 1$ the $r$-th osculating subspace $V^r_{\Gamma,p}$ is generated by 
$\Big\{\delta \big([\lambda^{\textrm{-}(2l\textrm{-}1)}]\big), l = 1,..,r\Big\}$. On the other hand, $\pi$ being separable,  the 
 tangent to $\iota_\pi(X)\subset Jac\,\Gamma$ at $0$ is equal to $\pi^*\big(H^1(X,O_X)\big)$, hence, generated by $\delta\big([\pi^*(z^{\textrm{-}1})]\big)$. 

In other words, the \textit{osculating order} $d$ is the smallest positive integer such that $\delta\big([\pi^*(z^{\textrm{-}1})]\big)$ is a linear combination $\sum^{d}_{l=1}a_l\delta\big([\lambda^{\textrm{-}(2l\textrm{-}1)}]\big)$, with $a_d\neq 0$. Or equivalently, thanks to the Mittag-Leffler Theorem, the smallest for which there exists a morphism $\kappa: \Gamma \to \mathbb{P}^1$, with polar parts equal to 
 $\pi^*(z^{\,\textrm{-}1})\,\textrm{-}\sum^{d}_{l=1}a_l\lambda^{\textrm{-}(2l\textrm{-}1)}$. The latter conditions on $\kappa$ 
 are equivalent to \textbf{2.6.}(1) \& (2). Moreover, up to replacing $\kappa$ by 
 $\frac{1}{2}\big(\kappa \,$-$\, \tau^*_\Gamma (\kappa)\big)$, we can assume $\kappa$ is 
 $\tau_\Gamma$-anti-invariant. The difference of two such functions should be 
 $\tau_\Gamma$-anti-invariant, while having a unique pole at $p$, of order strictly smaller than $2d\,$-$1\leq  2g\,$-$1$, where $g$ denotes the arithmetic genus of $\Gamma$.
 Hence the difference is identically zero, implying the uniqueness of such a morphism $\kappa$. $\blacksquare$\\
 
\textbf{ Definition 2.7.}
\begin {enumerate}
\item
\emph{The pair of marked projections $(\pi, \kappa)$, satisfying} \textbf{2.6.}(1),(2)\&(3), \emph{will be called  a hyperelliptic d-osculating pair, and $\kappa$ the hyperelliptic $d$-osculating
function associated to $\pi $.}
\item
\emph{If the latter $\pi :(\Gamma ,p)\rightarrow (X,q)$ does not dominate any other hyperelliptic $d$-osculating cover, we will call it
minimal-hyperelliptic $d$-osculating cover}.
\end{enumerate}

\textbf{Corollary 2.8.}\\
  \emph{Let $\pi: (\Gamma,p) \to (X,q) $ and $\pi': (\Gamma',p) \to (X,q) $ be two hyperelliptic covers of osculating orders, $d$ and $d'$ respectively, such that $\pi$ dominates $\pi'$. Then $d \leq d'$.}\\
  
\emph{Proof.} Let $\kappa'$ be the  \emph{hyperelliptic d-osculating} function associated to $\pi'$, and $ j:(\Gamma,p) \to (\Gamma',p')$ the birational morphism such that $\pi= \pi' \circ j$. Then, the poles of $ \kappa' \circ j: \Gamma \to \mathbb{P}^1$ lie along $\pi^{-1}(q)$, while  $\kappa'\circ j +\pi^*(z^{-1})=\big(\kappa'+{\pi'}^{\,*}(z^{-1})\big) \circ j $ has a pole of order $2d'\,$-$\,1$ and no other pole along $\pi^{-1}(q)$. It follows (along the same lines of proof as in \textbf{2.6.}) that the tangent to $\iota_\pi(X)$ must be contained in $V^{d'}_{\Gamma,p}$. Hence, the minimality of $d$ implies $d\leq d'$.$\blacksquare$

\section{The algebraic surface set up}

 \textbf{ 3.1.} We will construct hereafter the ruled surface $\pi_S: S \to X$ and its blowing-up $e: S^\bot \to S$, both naturally equipped with involutions $\tau : S \to S$ and $\tau^\bot: S^\bot \to S^\bot$, as well as a degree-$2$ projection $S^\bot \stackrel{\varphi}{\to} \widetilde S$ to a known anticanonical rational surface . We will then prove that any \emph{hyperelliptic $d$-osculating cover}
  $\pi: (\Gamma,p) \rightarrow (X,q)$  factors uniquely through $\pi_{S^\bot}: \pi_S \circ e: S^\bot \to X$ and projects, via $S^\bot \stackrel{\varphi}{\to} \widetilde S$, onto an irreducible rational curve. Moreover, we will prove that $\pi$ dominates a unique \emph{hyperelliptic $d$-osculating cover} (\textbf{3.9.}). \\

  \textbf{ Definition 3.2.}
  
\indent (1) \emph{Fix an odd meromorphic function $\zeta:X \to \mathbb{P}^1$, with divisor of zeroes and poles equal to $(\zeta)=q+\omega_1\,$-$\,\omega_2\,$-$\,\omega_3$, and consider the open affine subsets $U_o:=X\setminus \{q\}$ and $U_1:=X\setminus \{\omega_1\}$. We let $\pi_S: S \to X$ denote the ruled surface obtained by identifying $\,\mathbb{P}^1\times U_o\,$ with  $\,\mathbb{P}^1\times U_1\,$, over $X \setminus \{q,\omega_1\}$}:\\

 $\forall q' \neq q,\omega_1,\quad (T_o\,,q') \in  \mathbb{P}^1\times U_o\,\,\,$ \emph{is identified with} $\,\,\,(T_1+ \frac{1}{\zeta(q')}\,,q') \in  \mathbb{P}^1\times U_1$.\\

\emph{In other words, we glue the fibers of $\,\mathbb{P}^1\times U_0\,$ and $\,\mathbb{P}^1\times U_1\,$, over any $q' \neq q,\omega_1$, by means of a translation. In particular the constant sections $q' \in U_k \mapsto (\infty,q')\in \mathbb{P}^1\times U_k \,$ ($k=0,1$), get glued together, defining a particular one denoted by $C_o \subset S$}. \\

\indent (2)  \emph{The involutions $\,\mathbb{P}^1\times U_k \to\,\mathbb{P}^1\times U_k,\quad (T_k,q') \mapsto \big(\,$-$\,T_k,[\,$-$\,1](q')\big)\,\,\,(k=0,1)$, get glued under the above identification and define an involution $\tau: S \to S$, such that $\pi_S \circ \tau=[\,$-$\,1] \circ \pi_S$}. \emph{In particular, $\tau$ has two fixed points over each half-period $\omega_i$: one in $C_o$, denoted by $s_i$, and the other one denoted by $r_i$ ($i=0,..,3$). It can also be checked that translating along the fibers of $\,\mathbb{K}\times U_k$ by any scalar $a \in \mathbb{K}\,(k=0,1)$, extends to an automorphism $t_a:S \to S$, leaving fixed $C_o$ and such that $\pi_S \circ t_a =\pi_S$}.\\
 
 \indent (3) \emph{Whenever $\textbf{p}\geq 3$, we choose $\zeta$ }(\textbf{3.2.}(1)) \emph{as a local parameter of $X$ centered at $q$ , and consider the unique meromorphic function $f_{\textbf{p}}:X \to \mathbb{P}^1$, having a local development $f_{\textbf{p}}=\frac{1}{\zeta^{\textbf{p}}}+\frac{c}{\zeta}+O(\zeta)$, for some $c\in \mathbb{K}$. We denote  $C_{\textbf{p}}\subset S$ the curve defined over $\,\mathbb{P}^1\times U_o$ by the equation $\,T_o^{\textbf{p}} +cT_o+f_{\textbf{p}}=0\,$, and over $\,\mathbb{P}^1\times U_1\,$ by the equation 
$\,T_1^{\textbf{p}} +cT_1+f_{\textbf{p}}\,$-$\,\frac{1}{\zeta^{\textbf{p}}}\,$-$\,\frac{c}{\zeta}=0$}.\\

\textbf{Proposition 3.3.}\\
 \emph{The ruled surface $S \to X$ has a unique section of self-intersection $0$, namely $C_o$, and its canonical divisor is equal to -$2C_o$. In particular, $S \to X$ is isomorphic to $  \mathbb {P}(E)\to X $, the ruled surface associated to the unique indecomposable rank-$2$, degree-$0$ vector bundle over $X$}(cf. \cite{H}\S V.2, \cite{T-V}\S3.1.).\\
 
 \textbf{Proof.}  The meromorphic differentials $dT_o$ and $dT_1$ get also glued together, implying that $K_S$, the canonical divisor of $S$ is represented by -$2C_o$. Any section of $\pi_S: S \to X$, other than $C_o$, is given by two non-constant morphisms $f_i:U_i \to \mathbb{P}^1$ ($i=1,2$), such that $ \quad f_o=f_1\,$-$\,\frac{1}{\zeta}\quad$ outside $\{q,\omega_1\}$. 
A straightforward calculation shows that a section as above intersects $C_o$, while having self-intersection number greater or equal to $2$. It follows from the general Theory of Ruled Surfaces (cf. \cite{H}\S V.2) that $C_o$ must be the unique section with zero self-intersection. Hence, the ruled surface $\pi_S:S \to X$ defined above, is isomorphic to the projectivization of the unique indecomposable rank-$2$, degree-$0$ vector bundle over $X$(cf. \cite{H}\S V.2). $\blacksquare$\\

\textbf{Definition 3.4.}(cf. \cite{T-V}\S4.1.)\\
\emph{ Let $e :  S^\bot \to S $ denote hereafter the monoidal transformation of $S$ at
$\{s_i, r_i,i=0,..,3\}$, the eight fixed points of $ \tau$, and $\tau
^\bot : \ S^\bot \to S^\bot$ its lift to an involution fixing the
corresponding exceptional divisors $\big\{s_i^\bot:=e^{-1}(s_i), r_i^\bot:=e^{-1}(r_i), i=0,..,3\big\}$.
Taking the quotient of $S^\bot$ with respect to $\tau
^\bot$, we obtain a degree-$2$ projection $\varphi: S^\bot \to \widetilde S $, onto a smooth rational surface $\widetilde S $, ramified along the exceptional curves $\{s_i^\bot, r_i^\bot, i=0,..,3\}$. }\\
 
\textbf{Lemma 3.5.}\\
 \emph{ Whenever $\textbf{p}\geq 3$, the curve $C_{\textbf{p}}$} (\textbf{3.2.}(3)) \textit{is irreducible and linearly equivalent to $\textbf{p}C_o$. Moreover, any irreducible curve numerically equivalent to a multiple of $\,C_o$, is either $C_o$ itself or a translate of $C_{\textbf{p}}$. In particular  $C_{\textbf{p}}$ and $\textbf{p}C_o$ generate the complete linear system  $\big|\textbf{p}C_o\big|$, and $S$ is an elliptic surface.}\\
 
 \indent \textbf{Proof.} The curve $C_{\textbf{\emph{p}}}$ is $\tau$-invariant, does not intersect the section $C_o$ and projects onto $X$ with degree $\textbf{\emph{p}}$. Hence, $C_{\textbf{\emph{p}}}$ is linearly equivalent to $\textbf{\emph{p}}C_o$ and has multiplicity one at $r_o\in S$. In order to prove its irreducibility, we may assume $C_{\textbf{\emph{p}}} \to X$ is separable, or equivalently, that $c\neq 0$ in \textbf{3.2.}(3). Otherwise $C_{\textbf{\emph{p}}} \to X$ would be purely inseparable and $C_{\textbf{\emph{p}}}$ isomorphic to $X$. The curve $C_{\textbf{\emph{p}}}$ is then smooth and transverse to the fiber $S_o:= \pi_S^{-1}(q)$, and their intersection number at $r_o \in S_o\cap C_{\textbf{\emph{p}}}$ is equal to 1. Let $C'$ denote the unique irreducible $\tau$-invariant component of $C_{\textbf{\emph{p}}}$ going through $r_o$, and suppose that $C' \neq C_{\textbf{\emph{p}}}$. Then $C'$ has zero self-intersection and the projection  $C' \to X$ has odd degree $\textbf{\emph{p'}}$, for some $1<\textbf{\emph{p'}}<\textbf{\emph{p}}$. Otherwise (i.e.: if $\textbf{\emph{p'}}=1$), $C'$ would give another section of $\pi_S$ having zero self-intersection. Contradiction!  Its complement, say $C'':= C_{\textbf{\emph{p}}}\setminus C'$, is a smooth, effective divisor linearly equivalent to $(\textbf{\emph{p}}\,$-$\,\textbf{\emph{p'}})C_o$. Translating $C'$ by an appropiate automorphism $t_a$ (\textbf{3.2.}(2)), we may assume  that $t_a(C')$ intersects $C''$, hence $t_a(C') \subset C"$ because their intersection number is equal to $0$. It follows that any irreducible component of $C_{\textbf{\emph{p}}}$ is a translate of $C'$, forcing the prime number $\textbf{\emph{p}}$ to be a multiple of $p'>1$. Therefore, $\textbf{\emph{p}}=\textbf{\emph{p'}}$ and $ C_{\textbf{\emph{p}}}=C'$ is irreducible as asserted. 
Consider at last, any other irreducible curve, say $C$, linearly equivalent to $mC_o$ for some $m>1$. It has zero intersection number with $C_{\textbf{\emph{p}}}$ and must intersect some translate of $C_{\textbf{\emph{p}}}$, implying that they coincide. In particular $m=\textbf{\emph{p}}$ and any element of $\big|\textbf{\emph{p}}C_o\big|$, other than $\textbf{\emph{p}}C_o$, is a translate of $C_{\textbf{\emph{p}}}$.$\blacksquare$\\

  The Lemma and Propositions hereafter, proved in \cite{T.1}\S2.3.,\S2.4.,\& \S2.5., will be 
instrumental in constructing the equivariant factorization $\iota^\bot : \Gamma \to S^\bot$ (\textbf{3.2.}).\\ 

  \textbf{Lemma 3.6.} 

\noindent \textit{There exists a unique, $\tau$-anti-invariant, rational morphism $\kappa_s : S \to \mathbb{P}^1$,  
with poles over $C_o$}+\,\textit{$\pi_S^{\textrm{-}1}(q)$, such that over a suitable neighborhood 
$U$ of $q\in X$, the divisor of poles of $\kappa_s $}+\textit{\,$ \pi_S^*(z^{\textrm{-}1})$ is reduced and equal to
 $C_o\cap \pi_S^{\textrm{-}1}(U)$}. \\

  \textbf{Proposition 3.7.} 

\noindent \textit{For any hyperelliptic cover $\pi: (\Gamma, p) \rightarrow (X, q)$, the existence of the unique hyperelliptic $d$-osculating function $\kappa :\Gamma \rightarrow \mathbb{P}^1$ }(\textbf{2.7.}(1)) \textit{is equivalent to the existence of a unique morphism $\iota: \Gamma \to S$ such that $\iota \circ \tau_\Gamma= \tau \circ \iota$, $\,\pi = \pi_S \circ \, \iota\,$ and 
$\iota^*(C_o)=(2d\,$-$1)p$.}\\

 \textbf{Proposition 3.8.}\\
\noindent \textit{For any hyperelliptic $d$-osculating pair $(\pi,\kappa)$, the above morphism
$\iota: \Gamma \to S$ lifts to a unique equivariant morphism $\iota^\bot: \Gamma \to S^\bot$ 
(i.e.: $\tau^\bot \circ \iota^\bot = \iota^\bot \circ \tau_{\Gamma})$. 
In particular, $(\pi,\kappa)$ is the pullback of $(\pi_{S^\bot},\kappa_{s^\bot})=
(\pi_{S}\circ e,\kappa_{s}\circ e)$, and $\Gamma$ lifts to a $\tau^\bot$-invariant curve, $\Gamma^\bot:=\iota^\bot(\Gamma)\subset S^\bot$, which projects onto the rational irreducible curve $\widetilde \Gamma:=\varphi\big(\Gamma^\bot\big) \subset \widetilde S$. In particular, $2d\,$-$\,1= e^*(C_o)\cdot {\iota^\bot}_*(\Gamma)$.}
\begin{displaymath}
\xymatrix{
&  \Gamma^\bot \subset  S^\bot \ar[r]^\varphi \ar[d]_e \ar[ddr]^{\pi_{S^\bot}} & \widetilde \Gamma \subset   \widetilde S  \\
\Gamma \ar[r]|\iota \ar[ru]^{\iota^\bot} \ar[rrd]_\pi & \iota (\Gamma) \subset S \ar[rd]|{\pi_S} & \\
& & X
     }
\end{displaymath} 
 
 \textbf{Proof.}  The monoidal transformation $e : S^\bot \to S$, as well as $\iota : \Gamma \to S$, can be pushed down to the corresponding quotients, making up the following diagram:

 \begin{displaymath}
\xymatrix{
\Gamma  \ar[d]|{2:1} \ar[rd]^\iota &    &    S^\bot \ar[d]^\varphi \ar[ld]_e  \\
\Gamma/\tau_\Gamma  \ar[rd]^{\iota/} & S \ar[d]|{2:1}   &    \widetilde S \ar[ld]_{\widetilde e}  \\
                                                & S/\tau   &   }
\end{displaymath}

Moreover, since $\widetilde e : \widetilde S \to S/\tau$ is a birational morphism and $\Gamma/\tau_\Gamma$ is a smooth curve (in fact isomorphic to $\mathbb{P}^1\,$), we can lift  $ \iota/ :\Gamma/\tau_\Gamma \to S/\tau$  to  $\widetilde S$, obtaining a morphism    $\widetilde \iota:\Gamma   \to \widetilde \Gamma \subset \widetilde S$, fitting in the diagram:
\begin{displaymath}
\xymatrix{
&  \widetilde \Gamma \subset \widetilde S \ar[rd]^{\widetilde e} & \\
\Gamma \ar[ru]^{\widetilde \iota} \ar[rd]_\iota & & S/\tau \\
& S \ar[ru]|{2:1} &
  }
\end{displaymath}
 Recall now that $S^\bot$ is the fibre product of  $\,\widetilde e: \widetilde S \to S/\tau$  and  $ S \to S/\tau $ (cf. \cite{T-V}\S4.1.). Hence, $\,\iota\,$  and  $\,\widetilde \iota\,$  lift to a unique equivariant morphism $\iota^\bot : \Gamma \to \Gamma^\bot \subset S^\bot$, fitting in 

\begin{displaymath}
\xymatrix{
& &  \widetilde S \ar[rd]^{\widetilde e}\\
\Gamma \ar[rru]^{\widetilde \iota} \ar[rrd]_\iota \ar[r]|{\iota^\bot} & S^\bot \ar[ru]_\varphi \ar[rd]^e & & S/\tau\\
& &  S \ar[ru]|{2:1} }
\end{displaymath}

Furthermore, since $\widetilde \iota: \Gamma \rightarrow \widetilde S$ factors through $\Gamma \rightarrow \Gamma /\tau_{\Gamma}\cong \mathbb{P}^1$, its image  $\widetilde \Gamma:=\varphi\big(\iota^\bot(\Gamma)\big)=\widetilde {\iota} (\Gamma)  \subset \widetilde S$ is a rational irreducible curve as claimed.$\blacksquare$\\

\textbf{ Corollary 3.9.}\\
\noindent \textit{Any hyperelliptic $d$-osculating cover $\pi : (\Gamma ,p)\to (X,q)$ 
dominates a unique minimal-hyperelliptic d-osculating cover, with same image $ \Gamma^\bot \subset S^\bot$ as $\pi$.}\\

\textbf{Proof.} Let $\overline \pi : (\overline \Gamma ,\overline p)\to (X,q)$ be an arbitrary \textit{hyperelliptic $d$-osculating cover} dominated by $\pi : (\Gamma ,p)\to (X,q)$, $\psi: (\Gamma, p) \to (\overline \Gamma, \overline p)$ the corresponding birational morphism and $\overline \iota^\bot: \overline \Gamma \to S^\bot$ the 
factorization of $\overline \pi$ via $S^\bot$. The uniqueness of $\iota^\bot$ implies that 
$\iota^\bot = \overline{\iota}^\bot \circ \psi$. Hence, they have same image in $S^\bot$,  
$\iota^\bot(\Gamma) = \overline{\iota}^\bot(\overline \Gamma)= \Gamma^\bot$, and project onto the same curve  
$\widetilde \Gamma \subset \widetilde S$. Furthermore, $\psi$ and $\overline{\iota}^\bot$ being equivariant 
morphisms, we can push down $\psi: \Gamma \to \overline \Gamma$ to an identity between their quotients, 
$\Gamma/\tau_\Gamma \cong \mathbb{P}^1 \stackrel{=}{\rightarrow} \mathbb{P}^1\cong \overline \Gamma /\tau_{\overline \Gamma}$,
 as well as $\overline \iota^\bot$ to a morphism $\widetilde \iota : \mathbb{P}^1 \rightarrow \widetilde \Gamma$ (of same degree as $\overline{\iota}^\bot: \overline \Gamma \to \Gamma^\bot$), as shown hereafter:\\
 
\begin{displaymath}
\xymatrix{
 p \in \Gamma \ar@/_/[ddrr]|{2:1}   \ar[drr]_\psi \ar[drrr]^{\iota^\bot} \ar[rrrr]^\pi  &  & &  &  q \in X\\
                  &  &\overline p \in\overline \Gamma \ar[d]|{2:1}  \ar[r]_{\overline{\iota}^\bot} &  \Gamma^\bot \subset S^\bot \ar[d]^\varphi  \ar[ru]|{\pi_{S^\bot}} \\
&  & \mathbb{P}^1 \ar[r]^{\widetilde \iota}      & \widetilde \Gamma \subset \widetilde S   &\\
    }
\end{displaymath}

Taking the fiber product of  $\widetilde \iota: \mathbb{P}^1 \rightarrow \widetilde \Gamma$  and  $\varphi: \Gamma^\bot
\rightarrow \widetilde \Gamma$,  say  $\Gamma^\star$,  we then factorize $\overline \iota^\bot$ in the above diagram, through a birational morphism $\overline \Gamma  \to \Gamma^\star$  as follows:

\begin{displaymath}
\xymatrix{
    \overline p  \in \overline{\Gamma} \ar[dr] \ar[ddr]|{2:1} \ar[drr]^{\overline{\iota}^\bot} \ar[rrr]^{\overline \pi }& & & q \in X\\
 & p^\star \in \Gamma^\star \ar[d]|{2:1}    \ar[r]_{\iota^{\star \bot}}  & \Gamma^\bot \subset S^\bot \ar[d]^\varphi \ar[ru]|{\pi_{S^\bot}}  \\
               &   \mathbb{P}^1  \ar[r]^{\widetilde \iota}  & \widetilde \Gamma \subset \widetilde S  }
\end{displaymath} 

\noindent where $p^\star \in \Gamma^\star$ is the image of
$\overline p \in \overline \Gamma$. Furthermore, since $\overline p$ is smooth and the unique pre-image of $p^\star$, we deduce that the latter morphism factorizes via the desingularization of $\Gamma^\star$ at the unibranch point $p^\star$. We will therefore assume till the end of the proof, that $\Gamma^\star$ is indeed smooth at $p^\star$. On the other hand, the degree-$2$ projection ($\,\overline \Gamma \to \mathbb{P}^1$ is ramified at $\overline p$, hence) $\Gamma^\star \to \mathbb{P}^1$ is ramified at $p^\star$. Then, applying  \textbf{3.8.} one immediately checks that the natural projection $\pi^\star:= \pi_{S^\bot} \circ {\iota^\star}^\bot :(\Gamma^\star, p^\star) \rightarrow (X, q)$ is a \textit{hyperelliptic $d $-osculating cover}, dominated by $\overline \pi$ (and $\pi$ as well). Thus, the latter $\pi^\star$ is the unique \textit{minimal-hyperelliptic $d$-osculating cover} dominated by $\pi$.$\blacksquare$\\

 \textbf{ Remark 3.10.}
  
 \noindent The \textit{minimal-hyperelliptic $d$-osculating cover} $\pi^\star$, explicitely constructed in the proof of \textbf{3.9.}, can not be recovered from $\widetilde \Gamma:= \varphi(\Gamma^\bot)$, unless $m:= deg\,(\iota^\bot: \Gamma \rightarrow \Gamma^\bot)$ is equal to $1$. There exists indeed a $(m\,$-$\,1)$-dimensional family of (non-isomorphic) \textit{minimal-hyperelliptic $d$-osculating covers}, with same image $\widetilde \Gamma \subset \widetilde S$, as shown hereafter. We will actually start in \textbf{3.11.} from a \textit{minimal-hyperelliptic $d$-osculating cover} $\pi$ (i.e.: identifying $\,\Gamma\,$ with $\,\Gamma^\star\,$), and give its complete factorization, in terms of the rational curve $\widetilde \Gamma \subset \widetilde S$.\\

\textbf{ Corollary 3.11.}\\
\noindent \textit{ Let $\pi : (\Gamma , p ) \rightarrow (X, q)$ be a minimal-hyperelliptic $d$-osculating cover, equipped }(\textbf{3.8.})\textit{ with $\iota^\bot: \Gamma  \rightarrow  \Gamma^\bot $, its equivariant factorization through $S^\bot$, as well as $\,\mathbb{P}^1 \stackrel{j}{\rightarrow} \widetilde \Gamma$, the desingularization  of the rational irreducible curve $\widetilde \Gamma: =\varphi(\Gamma^\bot)$. Then, there exist unique marked morphisms $\psi :(\Gamma,p) \to(\Gamma^\flat,p^\flat)$, $\pi^\flat:(\Gamma^\flat,p^\flat)\to(X,q)$ and $\iota^{\flat \bot}:(\Gamma^\flat, p^\flat) \to \big(\Gamma^\bot, \iota^\bot(p)\big) $, such that (see the diagrams below):}\\

(1) \textit{$\pi$ and $\iota^\bot$ factor as $\pi^\flat \circ \psi$ and $ \iota^{\flat \bot}\circ \psi$,  respectively};\\

(2) $deg(\psi)=m:=deg(\iota^\bot)$\textit{, and $\psi^{\textrm{-}1}(p^\flat)= \{p\}$};\\

(3) \textit{$\pi^\flat$ is a minimal-hyperelliptic $d^\flat$-osculating cover, where $2d\,$-$\,1=m(2d^\flat$-$1)$}; \\

(4) \textit{there exist a polynomial morphism}  $\,R : (\mathbb{P}^1,\infty)  \stackrel{m:1}{\longrightarrow} (\mathbb{P}^1,\infty)$ \textit{and a degree-$2$ projection} $(\Gamma^\flat,p^\flat) \stackrel{f^\flat}{\to} (\mathbb{P}^1,\infty)$, \textit{such that $\Gamma $ is the fiber product of $\,R $ with }$f^\flat$ ;\\

(5) \textit{the arithmetic geni of $\,\Gamma $ and $\,\Gamma^\flat$, say $g $ and $g^\flat$,  satisfy $\,2g \,$}+$\,1 = m(2g^\flat$+$\,1)$.\\

(6) \textit{$\Gamma $ is isomorphic to $\Gamma^\bot$, if and only if, $m=1$ and $\widetilde \Gamma$ is isomorphic to $\mathbb{P}^1$.}\\

\textit{
Furthermore, the moduli space of degree-$n$ minimal-hyperelliptic $d$-osculating covers, having same image $\widetilde \Gamma \subset \widetilde S$ as $\pi $, is birational to a $(m$-$1)$-dimensional linear space}.\\

\textbf{Proof.} 
(1)-(2)-(3) Let $\Gamma^\flat$ denote the fiber product of $\Gamma^\bot \stackrel{\varphi}{\to}\widetilde \Gamma$ and $\mathbb{P}^1 \stackrel{j}{\rightarrow} \widetilde \Gamma$, equipped with the corresponding birational morphism $\Gamma^\flat \stackrel{\iota^{\flat \bot}}{\to} \Gamma^\bot$ and degree-$2$ cover $\Gamma^\flat \stackrel{f^\flat}{\to} \mathbb{P}^1$. The equivariant morphism $\iota^\bot$ can be pushed down, as in \textbf{3.9.}, to $\mathbb{P}^1 \stackrel{\widetilde \iota}{\to} \widetilde \Gamma$ and factors through $j$, say $\widetilde \iota= j \circ R$. Moreover, the latter morphisms satisfy $\varphi \circ \iota^\bot=\widetilde \iota= j \circ R$, implying the factorization through the fiber product $\Gamma^\flat$. In other words, there exists a degree-$m$ equivariant morphism $\Gamma \stackrel{\psi}{\to} \Gamma^\flat$ (i.e.: $\psi \circ  \,\tau_{\Gamma}=\tau_{\Gamma^\flat} \circ \,\psi$), such that $\iota^{ \bot}=   \iota^{\flat \bot} \circ \psi$, and with maximal ramification index at $p\in \Gamma$ (i.e.:  $\psi^{-1}(p^\flat)=\{p\}$, the fiber of $\iota^\bot$ over $\iota^\bot(p)$).
  In particular $\Gamma^\flat$ is unibranch at $p^\flat$, and up to replacing $(\Gamma^\flat, p^\flat)$ by its desingularization at $p^\flat$, we can assume $\pi^\flat:=\pi_{S^\bot} \circ \iota^{\flat \bot}: (\Gamma^\flat, p^\flat) \to(X,q)$ is a \emph{hyperelliptic cover}.
This construction is sketched in the diagrams below: 

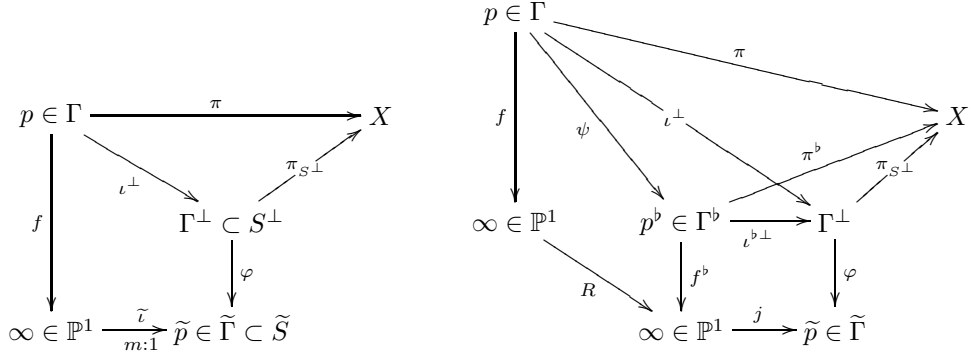
\begin{displaymath}
\xymatrix{& & & p \in \Gamma \ar[ddrr]|{\iota^\bot} \ar[dd]_{f} \ar[drrr]^\pi \ar[ddr]_{\psi} & & & \\
 p \in \Gamma \ar[dd]_{f} \ar[dr]_{\iota^\bot} \ar[rr]^\pi & & X & & &  & X\\ 
    &  \Gamma^\bot \subset S^\bot \ar[d]^\varphi \ar[ur]|{\pi_{S^\bot}} & &\infty \in \mathbb{P}^1 \ar[dr]_R & p^\flat \in \Gamma^\flat  \ar[d]^{f^\flat} \ar[r]_{\iota^{\flat \bot}} \ar[urr]^{\pi^\flat}  &  \Gamma^\bot  \ar[d]^\varphi \ar[ur]|{\pi_{S^\bot}} &  \\
                      \infty \in  \mathbb{P}^1 \ar[r]^{\widetilde \iota}_{m:1} & \widetilde p \in \widetilde \Gamma \subset \widetilde S & & &  \infty \in \mathbb{P}^1 \ar[r]^j & \widetilde p \in \widetilde \Gamma  &  }
\end{displaymath}

According to \textbf{3.8.}, the osculating order of $\pi^\flat$ (\textbf{2.4.}(2)), say $d^\flat$, satisfies $2d^\flat\,$-$\,1 =e^*(C_o) \cdot {\iota^{\flat \bot}}_*(\Gamma^\flat)$, while $2d\,$-$\,1 =e^*(C_o) \cdot {\iota^{ \bot}}_*(\Gamma)$. On the other hand, the factorization $\iota^\bot=  \iota^{\flat \bot}\circ \psi$ gives ${\iota^{ \bot}}_*(\Gamma)={\iota^{\flat \bot}}_*\big(\psi_*(\Gamma)\big)={\iota^{\flat \bot}}_*(m\Gamma^\flat)$, and replacing in the former equality gives $2d\,$-$\,1 =m(2d^\flat\,$-$\,1)$. Moreover, the \emph{minimal-hyperelliptic $d^\flat$-osculating cover} dominated by $\pi^\flat$ (\textbf{3.9.}) has same image $\Gamma^\bot$ as $\pi^\flat$, hence, it must dominate the fiber product product of $\Gamma^\bot \stackrel{\varphi}{\to}\widetilde \Gamma$ and $\mathbb{P}^1 \stackrel{j}{\rightarrow} \widetilde \Gamma$, and $\Gamma^\flat$ as well. In other words, $\pi^\flat$ is \emph{minimal-hyperelliptic}. \medskip
 
(4) Recall that $(\Gamma^\flat, p^\flat) \stackrel{f^\flat}{\longrightarrow} (\mathbb{P}^1, \infty)$ is classically represented in affine coordinates, as the zero locus $\big\{ y^2 = P(x)\big\}$ projecting onto the first coordinate, for some degree-$(2g^\flat$+$1)$ polynomial $P(x)$, $p^\flat$ being identified with the smooth Weierstrass point added at infinity. On the other hand, $\mathbb{P}^1\stackrel{R}{\to} \mathbb{P}^1$, the pushed down of $\Gamma \stackrel{\psi}{\to} \Gamma^\flat$ defined above, has maximal ramification index at $f(p)\in \mathbb{P}^1$ (i.e.: $f(p)\in \mathbb{P}^1$ is the unique pre-image of $f^\flat(p^\flat)\in \mathbb{P}^1$). Therefore, up to identifying the latter points with $\infty \in \mathbb{P}^1$, we may say that  $(\mathbb{P}^1,\infty)\stackrel{R}{\to} (\mathbb{P}^1,\infty)$ is defined by a degree-$m$ polynomial $R(t)$.
  Taking the fiber product of $\Gamma^\flat \stackrel{f^\flat}{\longrightarrow} \mathbb{P}^1$ with $ \mathbb{P}^1 \stackrel{R}{\longrightarrow} \mathbb{P}^1$, amounts then to replacing  $x$ by $R(t)$, giving the affine equation  $\big\{ y^2 = P\big(R(t)\big)\big\}$, where the composed polynomial $P\big(R(t)\big)$ has odd degree equal to $(2g^\flat$+$1)m$. Hence, the latter fiber product is a hyperelliptic curve, say $\Gamma_R$, of arithmetic genus $g_R$ such that $2g_R$+$1=m(2g^\flat$+$1)$, equipped with a smooth Weierstrass point $p_R \in \Gamma_R$ and a marked projection $(\Gamma_R, p_R) \stackrel{m:1}{\longrightarrow} (\Gamma^\flat, p^\flat)$, fitting in the following diagram:

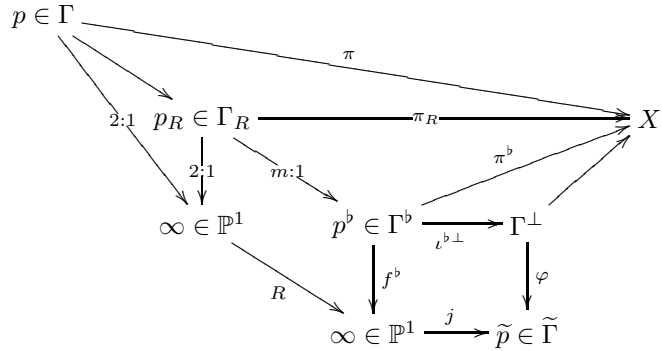
\begin{displaymath}
\xymatrix{  p \in \Gamma  \ar[drrrr]^\pi  \ar[dr] \ar[ddr]|{2:1} & & & \\
 & p_R \in \Gamma_R \ar[d]|{2:1} \ar[dr]|{m:1} \ar[rrr]|{\pi_R} & &  & X\\ 
    &  \infty \in \mathbb{P}^1 \ar[dr]_R & p^\flat \in \Gamma^\flat  \ar[d]^{f^\flat} \ar[r]_{\iota^{\flat \bot}} \ar[urr]^{\pi^\flat}  &  \Gamma^\bot  \ar[d]^\varphi \ar[ur] &  \\
                     &  &  \infty \in \mathbb{P}^1 \ar[r]^j & \widetilde p \in \widetilde \Gamma  &  }
\end{displaymath}

We can also check that $p_R \in \Gamma_R$ is the unique pre-image of $p^\flat \in \Gamma^\flat$, i.e.: the ramification index of $(\Gamma_R, p_R) \stackrel{m:1}{\longrightarrow} (\Gamma^\flat, p^\flat)$ at $p_R$ is equal to $m$. Hence, if $\kappa^\flat$ is the \textit{hyperelliptic $d^\flat$-osculating} function for $\pi^\flat$, its inverse image gives a \textit{hyperelliptic $d$-osculating} function for $\pi_R$. In other words, $\pi_R$ is a \textit{hyperelliptic $d$-osculating cover} dominated by the \textit{minimal-hyperelliptic $d$-osculating cover} $\pi$. Hence, they are isomorphic, implying that $\pi$ factors as $\pi^\flat \circ \psi$, $\,2g \,$+$\,1 = m(2g^\flat$+$\,1)$, and $\Gamma$ is the fiber product of $\mathbb{P}^1 \stackrel{R}{\longrightarrow} \mathbb{P}^1$ and $\Gamma^\flat \stackrel{f^\flat}{\longrightarrow} \mathbb{P}^1$, as claimed.\medskip

(5) It follows from the latter constructions that $\Gamma$ is isomorphic to $\Gamma^\bot$, if and only if $j:\mathbb{P}^1 \to \widetilde \Gamma$ is an isomorphism and $m=1$.\medskip

Consider at last, any other \textit{minimal-hyperelliptic $d$-osculating cover} having same image $\widetilde \Gamma \subset \widetilde S$. The latter must also factor through the above \textit{minimal-hyperelliptic $d^\flat$-osculating cover} $\pi^\flat$. We may replace then $R$ by any other degree-$m$ separable polynomial $P:\mathbb{P}^1 \to \mathbb{P}^1$, and take its fiber product with $\Gamma^\flat \stackrel{f^\flat}{\longrightarrow} \mathbb{P}^1$, to produce the general degree-$n$ \textit{minimal-hyperelliptic $d$-osculating cover} having image $\widetilde \Gamma$. Up to isomorphism, they are parameterized by a $(m$-$1)$-dimensional linear space. $\blacksquare$\\

 \section{ The hyperelliptic $d$-osculating covers as divisors of a surface}
\noindent \textbf{4.1.} The next step concerns studying the $\tau^\bot$-invariant irreducible curve $ \Gamma^\bot \subset S^\bot$, associated  in \textbf{3.} to any \emph{hyperelliptic cover} $\pi$. We calculate its linear equivalence class, in terms of the numerical invariants of $\pi$, and dress the basic relations between them. We also prove, whenever \textit{\textbf{p}}:=$char(\mathbb{K})\geq 3$, the supplementary bound $2g+1\leq \textit{\textbf{p}}(2d\,$-$\,1)$ \big(\textbf{4.4.}(1) \& (6)\big). We end up giving a numerical characterization for $\pi$ to be \emph{minimal-hyperelliptic} (\textbf{4.6.}). \\

 \textbf{ Definition 4.2.}\\
\noindent \textit{For any $i=0,..,3$, the intersection number between the divisors ${\iota^\bot}_*(\Gamma )$ 
and $r_i^\bot$ will be denoted by $\gamma _i$, and the corresponding vector $\gamma = (\gamma_i) \in \mathbb{N}^4$ called the type of $\pi$. Furthermore, for any $\mu=(\mu_i)\in \mathbb{N}^4$, $\mu ^{(1)}$ and $\mu ^{(2)}$ will denote, respectively:}\\

$\quad \quad \quad \quad \quad \mu ^{(1)} : = \sum^{3}_{i=0}\mu_i  \qquad$ and $\qquad \mu ^{(2)} : =  \sum^{3}_{i=0}\mu_i ^2.$\\
 
 \textbf{Lemma 4.3.} \\
 \noindent \textit{Let $(\Gamma, p) \stackrel{\pi}{\to} (X, q)$  be a degree-$n$ hyperelliptic $d$-osculating cover, of type $\gamma$ and ramification index $\rho$ at $p$. Consider its unique equivariant factorization through $S^\bot$, $\iota^\bot:\Gamma \rightarrow \Gamma^\bot$, and let $m$ denote its degree and $\iota := e \circ \iota^\bot$ its composition with the blowing up $S^\bot \stackrel{e}{\to} S$. Then :}

\begin {enumerate}

\item \textit{$ \iota_*(\Gamma )$ is equal to $m.\iota(\Gamma)$ and linearly
equivalent to $nC_o $}+\textit{$\,(2d$-$1)S_o$};

\item \textit{$\iota_*(\Gamma)$ is unibranch, and transverse to the fiber $S_o:=\pi_S^*(q)$, at $s_o= \iota(p)$};

\item \textit{$\rho$ is odd, bounded by $2d\,$-$\,1$ and equal to the multiplicity of $\iota_*(\Gamma)$ at $s_o$};

\item \textit{the degree $m$ divides $n$, $2d\,$-$\,1$ and $\rho$, as well as $\gamma_i$, for any } $i\in\{0,..,3\}$;

 \item $\gamma_o $+$ 1 \equiv \gamma_1 \equiv \gamma_2 \equiv \gamma_3 \equiv n(\textnormal{mod}.2)$;
  
  \item \textit{${\iota^\bot}_*(\Gamma)$ is linearly equivalent to 
 $e^*\big(nC_o$}+$(2d\,$-$1)S_o\big)\,$-$\,\rho \, s^\bot_o$ -$\,\sum^{3}_{i=0}\gamma_i\, r_i^\bot$.\\
 
\end{enumerate}

\textbf{Proof.} (1) Checking that $\iota_*(\Gamma)$ is numerically equivalent to $nC_o$+$(2d\,$-$1)S_o$ amounts 
to proving that the intersections numbers $\iota_*(\Gamma) \cdot S_o$ and $\iota_*(\Gamma)\cdot C_o$ are equal 
to $n$ and $2d\,$-$\,1$. The latter numbers are equal, respectively, to the degree of $\pi: \Gamma \rightarrow X$ 
and the degree of $\iota^*(C_o)= (2d\,$-$\,1)p$, hence the result. Finally, since $\iota_*(\Gamma)$ and $C_o$ only intersect at $s_o \in S_o$, we also obtain their linear equivalence.\\
  
 (2) $\&$ (3) Let $\kappa: \Gamma \to \mathbb{P}^1$ be the \textit{hyperelliptic $d$-osculating function} associated to $\pi$, uniquely 
characterized by properties \textbf{2.6.}(1),(2)\&(3), and $U \subset X$ a symmetric neighborhood of $q:=\pi(p)$. 
Recall that $\kappa + \pi^*(z^{\textrm{-}1})$ is $\tau_\Gamma$-anti-invariant and well defined over $\pi^{\textrm{-}1}(U)$, where it has a (unique) pole of order $2d\,$-$\,1$ at $p$. Studying its trace with respect to $\pi$ 
 we can deduce that $\rho$ must be odd and bounded by $2d\,$-$\,1$.
 
 On the other hand, let $\big(\iota_*(\Gamma), S_o\big)_{s_o}$
  and $\big(\iota_*(\Gamma), C_o\big)_{s_o}$ denote the intersection multiplicities at $s_o$, between 
  $\iota_*(\Gamma)$ and the curves $S_o$ and $C_o$. They are respectively equal, via the projection 
  formula for $\iota \,$, to $\rho$ and $2d\,$-$\,1$. At last, since $\iota_*(\Gamma)$ is unibranch at $s_o$
  and $\big(\iota_*(\Gamma), S_o\big)_{s_o} = \rho \leq 2d\,$-$\,1 = \big(\iota_*(\Gamma), C_o\big)_{s_o}$, 
  we immediately deduce that $\rho$ is the multiplicity of $\iota_*(\Gamma)$ at $s_o$ (and $S_o$ is transverse to $\iota_*(\Gamma)$ at $s_o$).
  
  (4) By definition of $m$, we clearly have $\iota_*(\Gamma) = m.\iota(\Gamma)$, while 
  $\{\rho, \gamma_i, i=0,..,3\}$ are the multiplicities of $\iota_*(\Gamma)$ at different points of $S$. 
  Hence, $m$ divides $n$ and $2d$-$1$, as well as all integers $\{\rho, \gamma_i, i = 0,..,3\}$.
 
  (5) For any $i= 0,..,3$, the strict transform of the fiber $S_i:= \pi_S^{-1}(\omega_i)$, by the monoidal transformation $e: S^\bot \to S$, is a $\tau^\bot$-invariant curve, equal to 
  $S^\bot_i:= e^*(S_i) \,$-$\, s_i^\bot \,$-$ \,r_i^\bot$, but also to $\varphi^*({\widetilde S}_i)$, where 
  ${\widetilde S}_i:= \varphi(S^\bot_i)$. Hence, the intersection number 
   ${\iota^\bot}_*(\Gamma)\cdot S_i^\bot$ is equal to the even integer \\

   ${\iota^\bot}_*(\Gamma)\cdot S_i^\bot = {\iota^\bot}_*(\Gamma)\cdot \varphi^*(\widetilde S_i)= 
   \varphi_*({\iota^\bot}_*\big(\Gamma)\big)\cdot \widetilde S_i = 2\widetilde \Gamma\cdot \widetilde S_i,$\\
   
   implying that $n= {\iota^\bot}_*(\Gamma)\cdot e^*(S_i)$ is congruent mod.$2$ to\\
    
 ${\iota^\bot}_*(\Gamma)\cdot S_i^\bot +\, {\iota^\bot}_*(\Gamma) \cdot (s_i^\bot +  \,r_i^\bot) \equiv {\iota^\bot}_*(\Gamma)\cdot (s_i^\bot +  \,r_i^\bot)($mod$.2).$\\
  
   We also know, by definition, that $\gamma_i:= {\iota^\bot}_*(\Gamma)\cdot r_i^\bot$, while ${\iota^\bot}_*(\Gamma)\cdot s_o^\bot= \rho$, the multiplicity of $\iota_*(\Gamma)$ at $s_o$, and  ${\iota^\bot}_*(\Gamma)\cdot s_i^\bot= 0$ if $i \neq 0$, because $s_i \notin \iota(\Gamma)$. Hence, $n$ is congruent mod.$2$, to $\rho $+$\, \gamma_o \equiv 1  $+$ \,\gamma_o\,($mod$.2)$, as well as to $\gamma_i$, if $i\neq 0$.\\
   
  (6) The Picard group $Pic(S^\bot)$ is the direct sum of $e^*(Pic(S))$ and the rank-$8$
  lattice generated by the exceptional curves $\{s_i^\bot, r_i^\bot, i= 0,..,3\}$. In particular, knowing that $\iota_*(\Gamma)$ is linearly equivalent to $nC_o $+$\, (2d\,$-$\,1)S_o$, and having already calculated  ${\iota^\bot}_*(\Gamma)\cdot s_i^\bot$ and ${\iota^\bot}_*(\Gamma)\cdot r_i^\bot$, for any $i= 0,..,3$, we can finally check that ${\iota^\bot}_*(\Gamma)$ is linearly equivalent to $e^*\big(nC_o $+$(2d\,$-$\,1)S_o\big)\, $-$\,\rho\,  s_o^\bot \,$-$\, \sum^{3}_{0}\gamma_i \, r_i^\bot.  \quad  \blacksquare$ \\

 \textbf{Theorem 4.4.} \\
\noindent \textit{Consider any hyperelliptic $d$-osculating cover $\pi:(\Gamma ,p)\to (X,q)$, of degree $n$, type $\gamma$, arithmetic genus $g$ and ramification index $\rho$ at $p$. Let $m$ denote the degree of its canonical equivariant factorization $\iota^\bot: \Gamma \rightarrow \Gamma^\bot\subset S^\bot$, and $\widetilde g$ the arithmetic genus of the rational irreducible curve $\,\widetilde \Gamma :=\varphi(\Gamma^\bot)$. Then, the numerical invariants $\{n, d, g, \widetilde g ,\rho, m, \gamma\}$ satisfy the following inequalities:}

\begin{enumerate}
\item $2g$+$1\leq \gamma ^{(1)}\quad$;\\

\item $4m^2\widetilde g = (2d$-$1)(2n$-$2m)+\,4m^2$-$\rho^2$-$\gamma^{(2)}\,$ and $\,\gamma^{(2)} \leq 2(2d\,$-$\,1)(n\,$-$\,m) $+$ 4\,m^2\,$-$\,\rho^2$;\\

\item $(2g$+$1)^2 \,\,\leq \,\,8(2d\,$-$\,1)(n\,$-$\,m) $+$ 13\,m^2$-$\,4\rho^2 \,\,\leq \,\,8(2d\,$-$\,1)n $+$ (2d\,$-$\,1)^2\quad$; \\

\item \textit{$\rho=1$ implies $m=1$, as well as } $(2g$+$1)^2\,\, \leq \,\, 8(2d\,$-$1)(n\,$-$\,1)$+$\,9\quad$;\\

\item \textit{if \textbf{p}$\,\geq3,\, $ we must also have} $\,\gamma^{(1)} \leq \textit{\textbf{p}}(2d\,$-$\,1)\quad$.\\
\end{enumerate}

\textbf{Proof.} (1) For any $i=0,..,3$, the fiber of $\pi_{S^\bot}:= \pi_S \circ e: S^\bot \rightarrow X$ over the half-period $\omega_i$, decomposes as $s_i^\bot $+$ \,r_i^\bot $+$\, S_i^\bot$, where $S_i^\bot$ is a $\tau^\bot$-invariant divisor and  $s_i^\bot$ is disjoint with ${\iota^\bot}_*(\Gamma)$, if $i \neq 0$, while ${\iota^\bot}^*(s_i^\bot)= \rho \, p$, by \textbf{4.3.}(2). Hence, the divisor $R_i:= {\iota^\bot}^*(r_i^\bot)$ of $\Gamma$ is linearly equivalent to $R_i\equiv \pi^{-1}(\omega_i)\,$-$\,(n\,$-$\,\gamma_i)\,p$ (and also $2R_i \equiv 2\gamma_i \,p\,)$. Recalling at last, that $\sum^{3}_{j=1}\omega_j \equiv 3\, \omega_o$, and taking inverse image by $\pi$, we finally obtain that $\sum^3_{i=0}R_i \equiv \gamma^{(1)}\, p\,$. In other words, there exists a well defined meromorphic function, (i.e.: a morphism), from $\Gamma$ to $\mathbb{P}^1$, with a pole of (odd!) degree $\gamma^{(1)}$ at the Weierstrass point $p$. The latter can only happen (by the Riemann-Roch Theorem) if $2g$+$1 \leq \gamma^{(1)}$, as asserted. \\

(2) The curve $\Gamma^\bot$ is $\tau^\bot$-invariant and linearly equivalent \big(\textbf{4.3.}(4)\&(6)\big) to:\\

\quad \quad $\Gamma^\bot\sim \frac{1}{m} \Big(e^*\big(nC_o $+$ \,(2d\,$-$1)S_o\big)\,$-$\,\rho  s_o^\bot $-$ \sum ^{3}_{i=0}\gamma_i\,  r_i^\bot\Big)$.\\

Recall also that $ \widetilde g\geq 0$ and $\widetilde K$, the canonical divisor of $\widetilde{S}$, is linearly equivalent to $\varphi _{*}\big(e^{*}(\,$-$\,C_0)\big)$ (\cite{T-V}\S4.2.(3)). Applying the projection formula for $S^\perp \stackrel{\varphi}{\to}\widetilde S$, to $\Gamma^\bot= \varphi^*(\widetilde \Gamma)$, we obtain
$\quad  0 \leq \widetilde g= \frac{1}{4m^2}\big((2d\,$-$\,1)(2n\,$-$\,2m)+\,4m^2\,$-$\,\rho^2 \,$-$\,\gamma^{(2)}\big)$, implying $\,\,\gamma^{(2)} \leq (2d\,$-$1)(2n\,$-$2m)$+$\, 4m^2$-$\,\rho^2 \,\, $, as claimed.
  
\medskip (3) \& (4) We start remarking that, for any $j= 1,2,3$, $(\gamma_o\,$-$\,\gamma_j)$ is a non-zero multiple of $m$. Hence,
$\sum_{\substack{i<j}}(\gamma_i\,$-$\,\gamma_j)^2 \geq 3m^2$, and replacing in \textbf{4.4.}(1) we get:
$$(2g\textrm{+}1)^2 \leq (\gamma^{(1)})^2 = 4\gamma^{(2)} \,\textrm{-} \,\sum_{\substack{i<j}}(\gamma_i\,\textrm{-}\,\gamma_j)^2 \leq  4\gamma^{(2)}\,\textrm{-}\,\,3m^2.
$$
Taking into account \textbf{4.4.}(3), we obtain the inequality \textbf{4.4.}(4). At last, since $m$ divides $\rho$ (\textbf{4.3.}(4)), $\rho=1$ implies $m=1$. Replacing in \textbf{4.4.}(3) gives us \textbf{4.4.}(4).\\

(5) Finally, let us assume \textit{\textbf{p}}$\,\geq3$ and denote by $C^\bot_{\textit{\textbf{p}}} \subset S^\bot$ the unique $\tau^\bot$-invariant irreducible curve, linearly equivalent to $ e^*(\textit{\textbf{p}}C_o)\,$-$\,\sum_{i=0}^3 r_i^\bot$. In particular, it can not be equal to $\Gamma^\bot$, hence  $C^\bot_{\textit{\textbf{p}}} \cdot \Gamma^\bot= \textit{\textbf{p}}(2d\,$-$\,1)\,$-$\,\gamma^{(1)}$ must be non-negative.$\blacksquare$\\ 
	
 \textbf{Corollary 4.5.} \\ 
 \noindent \textit{Let $\pi: \Gamma \rightarrow X$ be a degree-$n$ separable projection of a hyperelliptic curve onto the elliptic curve $X$, and let $g$ denote its arithmetic genus. Then, there exists a smooth Weierstrass point $p\in \Gamma$ such that $\pi: (\Gamma, p\,) \rightarrow \big(X, \pi(p)\big)$ is a hyperelliptic $d$-osculating cover, non ramified at $p$, with $d$ satisfying: $ (2d\,$-$\,1)(2n\,$-$\,2) \geq g^2+g\,$-$\,2$}.
 
\medskip \textbf{Proof.} Consider the global desingularization morphism $\overline j: \overline {\Gamma} \rightarrow \Gamma$, composed, either with  $\pi$, or with the degree-$2$ cover $\Gamma \rightarrow \Gamma/\tau_\Gamma \cong \mathbb{P}^1$.
As a ramified cover of $X$ and $\mathbb{P}^1$, we deduce from the Hurwitz formula that $\overline {\Gamma}$ is a smooth  hyperelliptic curve of positive genus, say $\overline {g}$, with $2\overline {g}\,$+$2$ Weierstrass points, while $\overline {\pi}:= \pi \circ \overline j: \overline{\Gamma} \rightarrow  X$ has, at most, $2\overline {g}\,$-$\,2$ ramifications points.  We can choose, therefore, a Weierstrass point $\overline{p} \in \overline{\Gamma}$, at which $\overline {\pi}$ is not ramified. In particular, its image $p:= \overline j(\overline {\,p})\in \Gamma$ must be a unibranch point. On the other hand, since $\overline {\pi}$ is not ramified at $\overline {p}$ and factors through $\pi: \Gamma \rightarrow X$, we see that $\pi$ restricts to a local isomorphism between neighborhoods of $p\in \Gamma$ and $q: =\pi(p) \in X$: 
$$\overline \pi:\,\,\overline p \in \overline \Gamma \stackrel{\overline j}{\to} p  \in \Gamma \stackrel{ \pi}{\to}q\in X$$

Hence, $p \,\,$ is a smooth Weierstrass point of $\Gamma$, at which $\pi$ is not ramified, and $\pi :(\Gamma,p) \to (X,q)$ is a \emph{hyperelliptic $d$-osculating cover} (\textbf{2.4.}(2)), for some integer $d \leq g$. Applying \textbf{4.4.}(4), we obtain $ (2d\,$-$\,1)(2n\,$-$\,2) \geq (g$+$\,2)(g\,$-$\,1)$ as claimed. $\blacksquare$\\

 \textbf{Corollary 4.6.}\\
\textit{Let $\pi: (\Gamma, p) \rightarrow (X, q)$ be a hyperelliptic $d$-osculating cover of type $\gamma$ and arithmetic genus $g$. Then $2g$}+\textit{$1 \leq \gamma^{(1)}$, with equality if and only if $\,\pi$ is minimal-hyperelliptic}.\\

 \textbf{Proof.} Recall that $\pi$ dominates a unique \textit{minimal-hyperelliptic $d$-osculating} (\textbf{3.9.}),  say $\pi^\star$, factoring through the same curve $\Gamma^\bot  \subset S^\bot$. Therefore, $\pi^\star$ has same type $\gamma$ as $\pi$, but a bigger arithmetic genus, say $g^\star$, satisfying $2g+1\leq2g^\star+1\leq\gamma^{(1)}$ (\textbf{4.4.}(1)). Hence, it is certainly enough to assume $\pi$ is \textit{minimal-hyperelliptic} and prove that $2g$+$1 \geq \gamma^{(1)}$.
 
Recall also, that $\iota^\bot:\Gamma \to \Gamma^\bot $ has odd degree $m$ and factors through the cover $\pi^\flat :(\Gamma^\flat, p^\flat) \to (X,q)$, of type $\gamma^\flat$ and arithmetic genus $g^\flat$, such that $\gamma^{(1)} = m\gamma^{\flat(1)}$ and $2g$+$1=m(2g^\flat$+$1)$ \big(\textbf{3.11.} \& \textbf{4.3.}(4)\big). Hence $2g\,$+$1= m(2g^\flat$+$1)\leq m\gamma^{\flat(1)}= \gamma^{(1)}$, with equality if and only if $2g^\flat$+$1= \gamma^{\flat(1)}$. We have thus reduced the problem, from $\pi$ to the \textit{minimal-hyperelliptic} $\pi^\flat$. So let us suppose in the sequel that $m=1$, or in other words, that $(\Gamma, p)=(\Gamma^\flat, p^\flat)$.
Let $(\Gamma^\lozenge, p^\lozenge) $ denote the fiber product of the marked morphisms $\big(\Gamma^\bot, \iota^\bot(p)\big) \stackrel{\varphi}{\longrightarrow}  (\widetilde \Gamma, \widetilde p \,)$ and $(\mathbb{P}^1, \infty) \stackrel{j}{\longrightarrow} (\widetilde \Gamma, \widetilde p \,)$ \big(\textbf{3.11.}\big). The marked curve $(\Gamma, p)=(\Gamma^\flat, p^\flat)$, is in fact the desingularization of $\Gamma^\lozenge$ at its unibranch point $ p^\lozenge$ (\textbf{3.11.}), and fits in the following diagram: 

\begin{displaymath}
\xymatrix{
  p \in \Gamma \ar[ddr]_{\iota^\bot} \ar[dr]|{1:1} \ar[ddd]_{\pi } \ar[drr]|{2:1} &  & \\
   &  p^\lozenge \in \Gamma^\lozenge\ar[r]|{f^\lozenge} \ar[d]|{1:1} &   \infty \in \mathbb{P}^1 \ar[d]^j\\
   & \iota^\bot(p) \in \Gamma^\bot  \ar[dl]|{\pi_{S^\bot}} \ar[r]^\varphi  & \widetilde p \in \widetilde \Gamma\\
   q \in X & & }
\end{displaymath}

\noindent Let $\widetilde g, g^\bot, g^\lozenge$ and $g $ denote the arithmetic geni of $\widetilde \Gamma, \Gamma^\bot,  \Gamma^\lozenge$ and $\Gamma $, respectively. Knowing the numerical equivalence class of $ \Gamma^\bot$ we easily obtain \big(e.g.: \textbf{4.4.}(2)\big):\\

$\quad \quad \widetilde g=\frac{1}{4}\Big((2d\,$-$1)(2n\,$-$2)$+$ \,4\,$-$\,\rho^2 \,$-$\,\gamma^{(2)} \Big) \quad$ and $\quad g^\bot = 2\widetilde g$+$\frac{1}{2}(\rho\,$-$\,2$+$\,\gamma^{(1)})$.\\

We can then deduce $g^\lozenge$, arguing as follows (like in the proof of \cite{T-V}\S5.8.(2)): since $ \Gamma^\bot \stackrel{\varphi}{\longrightarrow} \widetilde \Gamma$ is a flat degree-$2$ morphism, and $\mathbb{P}^1$ has arithmetic genus $=0$, we must have the relation $\, g^\bot\,$-$\,g^\lozenge\, = 2(\widetilde g\,$-$\,0)=2\widetilde g$. Hence, $g^\lozenge = \frac{1}{2}(\rho\,$-$\,2$+$\,\gamma^{(1)})$. We might as well argue that the desingularization morphism $\mathbb{P}^1 \stackrel{j}{\to} \widetilde \Gamma$ is obtained by monoidal transformation $\widetilde S$ (i.e.: $j\,$ is the restriction of a finite sequence of monoidal transformations $\ \widetilde S' \stackrel{j}{\longrightarrow} \widetilde S\,$ such that the strict transform of $\widetilde \Gamma \subset \widetilde S\,$ is isomorphic to $\mathbb{P}^1 $), implying that $\Gamma^\lozenge$ is contained in the fiber product of $S^\bot \stackrel{\varphi}{\longrightarrow} \widetilde S$ and $\widetilde S' \stackrel{j}{\longrightarrow} \widetilde S$, for which we can calculate its canonical divisor. Applying the adjunction formula gives the above  value of  $g^\lozenge$.

At last, composing $ (\Gamma,p)  \stackrel {1:1}\to (\Gamma^\lozenge,  p^\lozenge)$ with $(\Gamma^\lozenge, p^\lozenge) \stackrel{f^\lozenge}{\longrightarrow} (\mathbb{P}^1, \infty)$, we get the degree $2$ cover $f: \Gamma  \stackrel{f}{\to} \mathbb{P}^1$, and a morphism $( f,\pi ):\Gamma \to \Gamma_{f,\pi}  \subset \mathbb{P}^1\times X$ as in \textbf{2.5.}, fitting in:

\begin{displaymath}
\xymatrix{ 
& & \infty \in  \mathbb{P}^1\\
 p  \in \Gamma  \ar[r]|{1:1} \ar[rru]^f \ar[rrd]_{\pi } & p^\lozenge \in \Gamma^\lozenge \ar[ru]_{f^\lozenge} \ar[r] & ( \infty,q) \in \Gamma_{f,\pi }  \ar[u]|{2:1} \ar[d]|{n:1} \ar@{^{(}->}[r]&  \mathbb{P}^1\times X\\
 & &  q \in X
  }
\end{displaymath}

We have shown in the proof of \textbf{2.5.}(3), that $\frac{1}{2}(\rho\,$-$\,1)$ consecutive monoidal transformations are necessary to desingularize $\Gamma_{f,\pi}$ at its unibranch point $(\infty,q)$, and each monoidal transformation lowers its arithmetic genus by $1$. On the other hand, since $(\Gamma , p )$ dominates $(\Gamma^\lozenge,  p^\lozenge)$ and is smooth over $(\infty,q)$, we easily deduce that $g^\lozenge\,\textrm{-}\,g \leq \frac{1}{2}(\rho\,\textrm{-}1)$. Hence $g^\lozenge\,\textrm{-}\frac{1}{2}(\rho\,\textrm{-}1)= \frac{1}{2}(\textrm{-}1\textrm{+}\, \gamma^{(1)}) \leq g$. $\blacksquare$

\section{ On hyperelliptic $d$-osculating covers of arbitrary high genus}

 $\textbf{5.1.}$ -  We will let $C_o^\bot$ and $C_{\textit{\textbf{p}}}^\bot$ denote, hereafter, the strict transforms of $C_o$ and $C_{\textit{\textbf{p}}}$ by  $e:S^\bot \to S$ and $\widetilde C_o:= \varphi(C_o^\bot)$. Recall that to any \emph{hyperelliptic cover} $\pi:(\Gamma,p) \rightarrow (X,q)$ we have uniquely associated a morphism $\iota^\bot :\Gamma \to \Gamma^\bot \subset S^\bot$, a rational irreducible curve $\widetilde\Gamma:= \varphi(\Gamma^\bot) \subset \widetilde S$ and a vector $(n,d,\rho,\gamma) \in {\mathbb{N}^*}^3 \times \mathbb{N}^4$, satisfying the following restrictions (\textbf{4.3.} \& \textbf{4.4.}) :
 \begin{enumerate}
\item $\rho$ is odd, bounded by $2d\,$-$\,1$, and $\gamma_o $+$ 1 \equiv \gamma_1 \equiv \gamma_2 \equiv \gamma_3 \equiv n($mod$. 2)$;
\item if \textit{\textbf{p}} $\geq3$, we must have $\gamma^{(1)} \leq \textit{\textbf{p}}(2d\,$-$\,1)$.\\

Furthermore,  $\pi$ can be canonically recovered from $\widetilde \Gamma:=\varphi(\Gamma^\perp)$ if, and only if,  $\Gamma$ is birational to $\Gamma^\bot$, in which case:\\
 
\item  $\widetilde \Gamma$ has arithmetic genus $\widetilde g :=\frac{1}{4}\big((2d\,$-$\,1)(2n\,$-$\,2) $+$\, 4\,\,$-$\,\rho^2\,$-$\,\gamma^{(2)}\big)\geq0$;

 \item  $\Gamma^\perp=\varphi^*(\widetilde \Gamma)\,$ is linearly equivalent to $e^*\big(nC_o + (2d\,$-$1)S_o\big)\,$-$\,\rho {s_o}^\perp \,$-$\sum^{3}_{i=0}{\gamma_i {r_i}^\perp}$;

 \item  $ \widetilde \Gamma$ intersects $\widetilde{s}_o:=\varphi({s_o}^\perp)$, at a unique unibranch point, with multiplicity $\rho$;
 \item $\Gamma^\bot$ and $\widetilde\Gamma$ intersect $ C_o^\bot$ and $ \widetilde C_o$, (at most) at $ p_o^\bot:=C_o^\bot\cap s_o^\bot$ and $ \varphi(p_o^\bot)$, respectively, with multiplicities $2d\,$-$\,1\,$-$\,\rho$ and $\frac{1}{2}(2d\,$-$\,1\,$-$\,\rho)$.

 \end{enumerate}

 \textbf{Definition 5.2.} \\
 \textit{For any $(n, d, \rho, \gamma) \in {\mathbb{N}^*}^3 \times \mathbb{N}^4$ satisfying} \textbf{5.1.}(1),(2)\&(3)\emph{, we let $\Lambda(n,d,\rho, \gamma)$ denote the unique element of $Pic(\widetilde S)$ such that $\varphi^*\big(\Lambda(n,d,\rho, \gamma)\big)$ is linearly equivalent to $ e^*\big(nC_o + (2d\,$-$1)S_o\big)\,$-$\,\rho {s_o}^\perp \,$-$\sum^{3}_{i=0}{\gamma_i {r_i}^\perp}$, and $MH_X(n,d,\rho, \gamma)$ denote the moduli space of degree-$n$ minimal-hyperelliptic $d$-osculating covers of type $\gamma$, ramification index $\rho$ at their marked point, and birational to their canonical images in $S^\perp$}.\\
 
 \textbf{Proposition 5.3.}\\
 \textit{Any $\pi \in MH_X(n,d,\rho,\gamma)$ can be canonically recovered  from $\widetilde \Gamma 
\subset \widetilde S$ }\big(\textbf{3.11.}(2)\big)\textit{. Conversely, any rational irreducible curve $\widetilde \Gamma \subset \widetilde S$ satisfying properties} \textbf{5.1.}(1)-(6)\textit{, gives rise to a unique element of $MH_X(n,d,\rho,\gamma)$. }\\

\textbf{Proof.} Given $\widetilde \Gamma \subset \widetilde S$ satisfying \textbf{5.1.}(1)-(6), we denote $\Gamma^\bot:=\varphi^*(\widetilde \Gamma )\subset S^\bot$ and consider the fiber product of $(\Gamma^\bot, p^\bot) \stackrel{\varphi}{\to} \big(\widetilde \Gamma,\varphi(p^\bot)\big) $ with the desingularization morphism $(\mathbb{P}^1,\infty) \stackrel{j}{\to}\big(\widetilde \Gamma,\varphi(p^\bot)\big) $, say $(\Gamma,p)$. Proceeding as in the proof of \textbf{3.11.}, for the construction of $\pi^\flat$, we can easily prove that the natural domination $(\Gamma,p)\to (\Gamma^\bot,p^\bot)$, composed with $\pi^\bot:(\Gamma^\bot,p^\bot)\to(X,q)$ is indeed the announced \emph{minimal-hyperelliptic $d$-osculating cover}.$\blacksquare$\\

Studying $MH_X(n,d,\rho, \gamma)$ for a general vector $(n,d,\rho, \gamma)$, is a difficult and elusive problem. 
We will henceforth restrict to the simpler case where $\rho=1$ and $\widetilde \Gamma$ is isomorphic to $\mathbb{P}^1$. In other words, we will focus on degree-$n$ \textit{minimal-hyperelliptic $d$-osculating covers} with $\rho=m=1$, and type $\gamma$ satisfying $\gamma^{(2)}=(2d\,$-$\,1)(2n\,$-$\,2)\,$+$\,3$ (as well as \textit{$\gamma^{(1)} \leq \textit{\textbf{p}}(2d\,$-$\,1)$, if \textbf{p}}$\geq3$). \\

\textbf{Proposition 5.4.} (\cite{T.1}\S3.4)\\
\textit{Any curve $\Gamma \subset S$ intersecting $C_o $ at a unique smooth point $p\in\Gamma$ is irreducible, unless \textit{\textbf{p}}$\,\geq3$ and $C_{\textit{\textbf{p}}}$ is a component of $\Gamma$.}\\

\textbf{Proposition 5.5.}\\
\textit{
 Let $\Gamma ^\bot \subset S^\bot$ be a curve with no irreducible component in $\{r_i^\bot, i=0,..,3\}$, and intersecting $C_o^\bot $ (at most) at a unique smooth point $p^\bot\in\Gamma^\bot$. 
Then, $\Gamma ^\bot$ is an irreducible curve, unless $\textit{\textbf{p}}\geq3$ and $C_{\textit{\textbf{p}}}^\bot$ is a component of $\,\Gamma^\bot$.}\\

\textbf{Proof.} The properties satisfied by $\Gamma^\bot$ assure us that $\Gamma : = e _*(\Gamma^\bot)$, its direct image by
$e:S^\bot \to S$, does not contain $C_o$, and that $\Gamma^\bot$ is the
strict transform of $\Gamma$. We can also check, that $\Gamma $ is smooth at $p:= e(p^\bot)$ and $\Gamma \cap C_o=\{p\}$. It follows, by \textbf{5.4.},
that ($\Gamma $, as well as its strict transform) $\Gamma^\bot$ is, either an irreducible curve, or $\textit{\textbf{p}}\geq3$ and $C_{\textit{\textbf{p}}}^\bot$ is a component of $\Gamma^\bot$.$\blacksquare$\\

 \textbf{Proposition 5.6.} (\cite{T-V}\S6.2. \& \cite{L})

 \noindent \textit{Any $\alpha  = (\alpha _i)\in \mathbb
N^4$ such that $\,\alpha^{(2)}=2a+1$ is odd (and $\alpha^{(1)}\leq \textit{\textbf{p}}$, whenever $\textit{\textbf{p}}\geq3$), gives rise to an exceptional curve of the first kind $\widetilde \Gamma_\alpha \subset \widetilde S$. More precisely, let $k\in \{0,1,2,3\}$ denote the index satisfying $\alpha _k+1 \equiv\alpha _j$\textnormal{(mod.2)}, for any
$j\neq k$, and $S_k:=\pi_S^{-1}(\omega_k)$, then $\widetilde \Gamma_\alpha\,$ is a $(\,$-$\,1)$-curve and $\varphi^*(\widetilde \Gamma_\alpha) \subset S^\bot$ is the unique $\tau^\bot$-invariant irreducible curve linearly equivalent to
$e^*(aC_o+S_k)\,$-$\,s_k^\bot$-$\sum _{i=0}^3\alpha _i r_i^\bot$}. \\

\textbf{Proof.} Let $\Lambda$ denote the unique numerical equivalence class of $\widetilde S$ satisfying $\varphi^*(\Lambda)=e^*(aC_o+S_k)\,$-$\,s_k^\bot$-$\sum _{i=0}^3\alpha _i r_i^\bot$. It has self-intersection $\Lambda \cdot \Lambda=\,$-$\,1$, and $\Lambda\cdot \widetilde K=\,$-$\,1$ as well, hence, $h^o\big(\widetilde S, O_{\widetilde S}(\Lambda)\big)\geq \chi\big(O_{\widetilde S}(\Lambda)\big)=1$, and there exists an effective divisor $\widetilde \Gamma \in \big|\Lambda\big|$. If $\textit{\textbf{p}}=0$, such a divisor $\widetilde \Gamma$ is known to be unique and irreducible (\cite{T-V}\S6.2.). Its proof takes in account that for any $m>1$ there is no irreducible curve in $S$, numerically equivalent to $mC_o$. However, when $\textit{\textbf{p}}\geq3$ the latter property fails, due to the existence of ${C_{\textit{\textbf{p}}}}\subset S$, implying that the intersection number $C_{\textit{\textbf{p}}}\cdot \Lambda=\textit{\textbf{p}}\,\,$-$\, \alpha^{(1)}$ must be non-negative. Conversely, if $\alpha^{(1)}\leq \textit{\textbf{p}}$, $\Lambda$ intersects non-negatively ${\widetilde C}_{\textit{\textbf{p}}}:=\varphi({C_{\textit{\textbf{p}}}}^\bot)$, (as well as all other (-1) and (-2)-curves in $\widetilde S$), and M.Lahyane's irreducibility criterion for (-1)-classes applies to $\Lambda$ (\cite{L}).
$\blacksquare$\\

\indent According to \textbf{5.6.}, any $\alpha \in \mathbb{N}^4$ such that $\alpha^{(2)}$ is odd (and $\alpha^{(1)}\leq \textit{\textbf{p}}\,$, if $\,\textit{\textbf{p}}\geq3$), gives rise to an
exceptional curve of the first kind $\widetilde \Gamma_\alpha \subset  \widetilde S$. Conversely, we have the \\

\textbf{Corollary 5.7.} \\   
\emph{Any irreducible curve in  $\widetilde S$, with negative self-intersection, is either equal to some  $\widetilde \Gamma_\alpha$ } as above (\textbf{5.6.}), to $\widetilde C_{\textit{\textbf{p}}}$ if $\textit{\textbf{p}}\geq3$, \emph{or belongs to the set $\Big\{\widetilde C_o, \widetilde s_i, \widetilde r_i, i=0,..,3 \Big\}$}. \\ 

\textbf{Proof.} The arithmetic genus of  an arbitrary irreducible curve $\widetilde \Gamma\subset \widetilde S$  is non-negative and equal to $\widetilde g:= 1\,$+$\,\frac{1}{2}\big(\widetilde \Gamma \cdot  \widetilde \Gamma\,$+$\,\widetilde \Gamma \cdot  \widetilde K\big) \geq0$, where $\widetilde K$ denotes the canonical divisor of $\widetilde S$.  In particular $\widetilde \Gamma \cdot  \widetilde \Gamma\,$+$\,\widetilde \Gamma \cdot  \widetilde K\geq\,$-$\,2$. Moreover, since  $\varphi^*(\widetilde K)=e^*($-$\,2C_o)$  (cf. \cite{T-V}) and $C_o$ is \textit{nef}, we immediately deduce that $\widetilde \Gamma \cdot  \widetilde K\leq0$. Hence, $\,\widetilde \Gamma\cdot \,\widetilde \Gamma<0$ implies, either $\widetilde \Gamma \cdot  \widetilde \Gamma=\,$-$\,2$ and $\widetilde \Gamma \cdot  \widetilde K=0$, or $\widetilde \Gamma \cdot  \widetilde \Gamma=\,$-$\,1=\widetilde \Gamma \cdot  \widetilde K$. It follows, in any case, that $\widetilde g=0$, hence $\widetilde \Gamma$ is isomorphic to $\mathbb{P}^1$. If $\widetilde \Gamma\cdot \widetilde \Gamma=\,$-$\,1=\widetilde \Gamma \cdot \widetilde K$, one can easily check, via the projection formulae for $ S^\bot \stackrel{\varphi}{\to} \widetilde S$ and $S^\bot \stackrel{e}{\to} S$, that $\Gamma^\bot:=\varphi^*(\widetilde \Gamma) $ is a $\tau^\bot$-invariant divisor in $S^\bot$ and its projection in $S$, $\Gamma:=e_*(\Gamma^\bot)$, satisfies:\\

\indent $\Gamma \cdot C_o= e_*(\Gamma^\bot)\cdot C_o=\Gamma^\bot \cdot e^*(C_o)=\,$-$\,\frac{1}{2}\Gamma^\bot\cdot e^*(\,$-$\,2C_o)=\,$-$\,\frac{1}{2}\Gamma^\bot\cdot \varphi^*(\widetilde K)= \,$-$\, \widetilde\Gamma\cdot \widetilde K =1$. \\

It immediately follows that $\Gamma$ (as well as $\Gamma^\bot$) is irreducible. Otherwise it would break as a sum of two divisors exchanged by $\tau:S \to S$, in which case the above intersection number $\Gamma \cdot C_o$ should have been even. In other words, $\Gamma$ is an irreducible $\tau$-invariant curve, intersecting $C_o$ at $s_k$, for a unique $k \in \{0,1,2,3\}$.
Hence, $\Gamma$ is linearly equivalent to $aC_o$+$S_k$, for some $a\in \mathbb{N}$. 

Recall also that  $\Gamma^\bot \cdot (C_o^\bot $+$\sum_{i=0}^3 s_i^\bot)=\Gamma^\bot \cdot e^*(C_o)=1$, and let $\alpha=(\alpha_i)$ denote the vector of intersection numbers  $( \Gamma^\bot \cdot r_i^\bot)$. Then, $\Gamma^\bot$ is linearly equivalent to $e^*(aC_o$+$S_k)\,$-$\,s_k^\bot\,$-$\,\sum_{i=0}^3 \alpha_ir_i^\bot$, and intersecting with the numerically equivalent curves $\big\{S_i^\bot:=e^*(S_i)\,$-$\,s_i^\bot\,$-$\,r_i^\bot, \,i=0,1,2,3\big\}$ one easily finds  out that $\alpha_k$+$1\equiv \alpha_i$(mod.2), for any $i\neq k$. Moreover, its self-intersection is equal to\\

\indent $2a\,$-$\,1\,$-$\,\alpha^{(2)}=\Gamma^\bot \cdot \Gamma^\bot =\varphi^*(\widetilde \Gamma) \cdot \Gamma^\bot =\widetilde\Gamma \cdot \varphi_*(\Gamma^\bot)=2\widetilde \Gamma\cdot \widetilde \Gamma= \,$-$\,2$. \\

 In other words, $2a+1=\alpha^{(2)}$ and $\widetilde \Gamma = \widetilde \Gamma_\alpha$ (\textbf{5.6.}).
 
At last, let us suppose that $\widetilde \Gamma \cdot\widetilde \Gamma=\,$-$\,2$ and $\widetilde \Gamma \cdot \widetilde K=0$, but $\widetilde \Gamma$ does not belong to $\big\{\widetilde s_i, \widetilde r_i, i=0,..,3 \big\}$. It then follows that $\Gamma^\bot:=\varphi^*(\widetilde \Gamma)$ is a $\tau^\bot$-invariant divisor of  $S^\bot$, of self-intersection $\Gamma^\bot\cdot\Gamma^\bot=\,$-$\,4$, equal to the strict transform of $\Gamma:=e(\Gamma^\bot) \subset S$. Therefore, it must be, either an irreducible degree-$2$ cover of $\,\widetilde \Gamma$, or break as the sum of two copies of $\widetilde \Gamma\simeq \mathbb{P}^1$, interchanged by $\tau^\bot$. In the latter case, $\Gamma^\bot$ should be the strict transform of  the divisor ${\pi_S}^{-1}(q'+[\,$-$\,1]q')$, for some $q' \in X$, in which case $\Gamma^\bot\cdot\Gamma^\bot\neq\,$-$\,4$. Hence, $\Gamma^\bot$ is indeed irreducible (and $\Gamma=e_*(\Gamma^\bot)$ as well). On the other hand,  recalling that $\varphi^*(\widetilde K)=e^*( \,$-$\,2C_o)$ and  $\varphi_*(\Gamma^\bot)=2\widetilde\Gamma$, we obtain\\

$\Gamma\cdot (\,$-$\,2C_o)=e_*(\Gamma^\bot) \cdot (\,$-$\,2C_o)=\Gamma^\bot \cdot e^*(\,$-$\,2C_o)=\Gamma^\bot \cdot \varphi^*(\widetilde K)=2\widetilde \Gamma\cdot \widetilde K=0$,\\

\noindent implying $\Gamma$ is numerically equivalent to a multiple of $C_o$ According to $\textbf{3.5.}$ this can only happen if $\Gamma=C_o$ and $\Gamma^\bot=C_o^\bot$, or \textit{\textbf{p}}$\geq3$, $\Gamma=C_{\textit{\textbf{p}}}$ and $\Gamma^\bot=C_{\textit{\textbf{p}}}^\bot$. $\blacksquare$\\ 

\textbf{Lemma 5.8.}\\
\emph{Let $ \Lambda:=\Lambda(n,d,1,\gamma)$ be as in} \textbf{5.2.}, $\widetilde \Gamma $ \textit{an arbitrary exceptional curve of the first kind on $\widetilde S$ and $\alpha \in \mathbb{N}^4$ the unique vector as in }\textbf{5.6.} \textit{such that $\widetilde \Gamma =\widetilde \Gamma_\alpha$ }(\textbf{5.7.}).\textit{ Then:}

\begin{displaymath}
4(2d\textrm{-}1)\widetilde \Gamma_\alpha \cdot \Lambda=\left\{\begin{array}{lll}
 \big(\gamma\textrm{-}\,(2d\,\textrm{-}\,1)\alpha\big)^{(2)}\,\,\textrm{-}\,\,\,(2d\,\textrm{-}\,1)^2\,\textrm{-}\,\,3 \quad, \quad if \quad\Gamma_\alpha \cdot \widetilde s_o=1\\ \\
\medskip\big(\gamma\,\textrm{-}(2d\,\textrm{-}\,1)\alpha\big)^{(2)}+2(2d\,\textrm{-}\,1)\,\textrm{-}\,(2d\,\textrm{-}\,1)^2\,\textrm{-}\,\,3 \quad otherwise.
\end{array} \right.
\end{displaymath}

\textbf{Proof.} Straightforward verification.$\blacksquare$\\ 

For $\Lambda(n,d,1,\gamma)$ to be \emph{nef}, we must have $\Lambda(n,d,1,\gamma)\cdot \widetilde \Gamma_\alpha \geq 0$, for any $\alpha$ as above. On the other hand, minimizing their value is tantamount (\textbf{5.8.}) to minimizing the norm of  $\gamma\,$-$\,(2d\,$-$\,1)\alpha$. In order to do it we make the following definitions.\\

\textbf{Definition 5.9.}

\begin{enumerate}
\item\emph{Given $(n,d,\gamma) \in \mathbb{N}^*\times\mathbb{N}^*\times\mathbb{N}^4$ satisfying} $\gamma_o +1 \equiv \gamma_j$(mod.$2), \forall j=1,2,3$\emph{, as well as} $\gamma^{(2)}=(2d\,$-$\,1)(2n\,$-$\,2)+3$\emph{, we let $\gamma = (2d\,$-$\,1)\mu+2\varepsilon$ be the  unique decomposition, with $\mu \in \mathbb{N}^4$ having same parity as $\gamma$, and $\varepsilon \in \mathbb{Z}^4$ such that} max$\{ |\varepsilon_i|\} \leq d\,$-$\,1$.  \emph{We will also assume, here and henceforth, that $\gamma^{(1)}=(2d\,$-$\,1)\mu^{(1)}+2\varepsilon^{(1)}\leq \textit{\textbf{p}}(2d\,$-$\,1)$, whenever} $\textit{\textbf{p}}\geq3$. \\

\item \emph{We define $^\natural  \mu=(^\natural  \mu_i) \in \mathbb{N}^4$ in order to have $(^\natural  \mu_i\,$-$\,\mu_i)\varepsilon_i=|\varepsilon_i|\,,\,\forall i=0,\cdots,3$}:
$$^\natural  \mu_i=\mu_i+1\quad\quad \textnormal{if} \quad \quad\varepsilon_i \geq0 \quad\quad \textnormal{or}\quad \quad ^\natural  \mu_i=\mu_i-1\quad\quad \textnormal{if} \quad \varepsilon_i <0$$

\item \emph{At last, we choose two indices $i_o\neq j_o$, where $|\varepsilon_i|$ attains its two maximal values, and  let $^\flat  \mu =(^\flat  \mu_i)\in \mathbb{N}^4$ be such that for all $ i\in \{0,1,2,3\}:$}
$$\quad ^\flat   \mu_i=^\natural\mu_i\quad\quad \textnormal{if} \quad \quad i \in \{i_o,j_o \} \quad\quad \textnormal{or}\quad \quad ^\flat  \mu_i=\mu_i\quad\quad \textnormal{if} \quad i \notin \{i_o,j_o \}$$

\end{enumerate}

\textbf{Remark 5.10.}

The vector $^\flat \mu $ may not be uniquely defined by \textbf{5.9.}(3). 
It  should also be clear that $\mu \equiv \gamma$(mod.$2$), and $4\varepsilon^{(2)} \equiv 3 \big($mod.($2d\,$-$\,1)\big)$. Conversely, we have the \\

\textbf{Proposition 5.11.}\\
\textit{Given any $n\in \mathbb{N}^*$  and $\gamma =(2d\textnormal{-}1)\mu+2\varepsilon$, with $\mu \in \mathbb{N}^4$ and $\varepsilon\in \mathbb{Z}^4$, such that: }
\begin{align}
 & \mu_o +1 \equiv \mu_1\equiv \mu_2 \equiv \mu_3(\textnormal{mod}.2)   \notag \\
& 4\varepsilon^{(2)} \equiv 3 \big(\textnormal{mod}.(2d-1)\big)    \quad \textnormal{and} \quad|\varepsilon_i|\leq d\,\textnormal{-}1, \quad i=0,\cdots, 3, \notag \\
& \gamma^{(2)}=(2d\textnormal{-}1)(2n\textnormal{-}2)+3, \notag \\
 & \textnormal{(as well as} \quad\gamma^{(1)}\leq \textit{\textbf{p}}(2d\,\textnormal{-}\,1\,), \quad \textnormal{if} \,\,\textit{\textbf{p}}\geq3 ), \notag
\end{align}

\noindent \emph{the minimal value of $\Lambda(n,d,1,\gamma)\cdot \widetilde \Gamma_\alpha$, taken amongst all $\alpha \in  \mathbb{N}^4$ with $\alpha^{(2)}$ odd, is attained at $\alpha$ equal, either to $\mu$, to $^\natural \mu$, or to $^\flat \mu$.}\\

\textbf{Corollary 5.12.}\\
\emph{The divisor $\Lambda(n,d,1,\gamma)$ is nef if and only if the vector $2\varepsilon=\gamma\,$-$\,(2d\,$-$\,1)\mu \in \mathbb{Z}^4$ \textnormal{(\textbf{5.9.})}, such that $4\varepsilon^{(2)}\equiv 3\big($mod$.(2d\,$-$\,1)\big)$ and} max$\{|\varepsilon_i|\}\leq d\,$-$\,1$, \emph{satisfies the supplementary conditions :}

\begin{enumerate}

\item $\varepsilon^{(2)}\geq d^2\,$-$\,d+1$;\\

\item $(2d\,$-$\,1)(^\natural \mu\,$-$\,\mu)\cdot \varepsilon=(2d\,$-$\,1)\big(\sum_{i=0}^{3}|\varepsilon_i|\big)\leq 3d^2\,$-$\,3d+\varepsilon^{(2)};$\\

\item $(2d\,$-$\,1)(^\flat \mu\,$-$\,\mu)\cdot \varepsilon = $ max$\big\{|\varepsilon_i|$+$|\varepsilon_j|,\forall i\neq j,\big\}\leq d^2\,$-$\,1+\varepsilon^{(2)}$.\\

\end{enumerate}

  As we shall see, given any $n,d \in \mathbb{N}^*$, there exist \emph{types} $\gamma =(2d\,$-$\,1)\mu$+$2\varepsilon \in \mathbb{N}^4$, such that $\gamma_o$+$1\equiv \gamma_1\equiv \gamma_2 \equiv \gamma_3$(mod.2) and $\gamma^{(2)}=(2n\,$-$\,2)(2d\,$-$\,1)$+$3$, for which $\Lambda:=\Lambda(n,d,1,\gamma)$ is, either \emph{nef} or not. We will actually construct in \textbf{5.13.} and \textbf{5.14.}, explicit examples where, either $\varepsilon$ satisfies \textbf{5.12.}(1),(2) \&(3), hence $\Lambda$ is \emph{nef}, or it does not satisfy \textbf{5.12.}(1), hence $\Lambda$ is not \emph{nef}. We actually conjecture that \textbf{5.13.} exhausts all \emph{types} such that $\gamma^{(2)}=(2d\,$-$\,1)(2n\,$-$\,2)\,$+$\,3\,$ and $\,\Lambda(n,d,1, \gamma)$ is \emph{nef}. \\

\textbf{Proposition 5.13}\\
\emph{Let us fix $\,d\geq2\,$, $\,k\in \{0,1,2,3\}\,$,  and $\mu\in \mathbb{N}^4$ such that $\mu_o+1\equiv
\mu_j \textnormal{(mod.2)}$ (for $j=1,2,3$). Pick any vector $
 \,2\varepsilon=(2\varepsilon_i)\in 2\mathbb{Z}^4\,$, satisfying} $(\forall i=0,\ldots ,3)\,$:

\begin{displaymath}
either \quad\vert 2\varepsilon_i\vert = (2d\textrm{-}2)(1\textrm{-}\delta
_{i,k})\,,\quad
or \quad\,\, \left\{\begin{array}{lll}
 |2\varepsilon_i| = d\,\,\textrm{-}\,(\textrm{-}1)^{\delta_{i,k}}\quad\quad \textit{if} \quad d \quad\textit{is odd}\quad,\\ \\
\medskip |2\varepsilon_i| = d\,\textrm{-}\,2 \delta_{i,k}\quad\quad \quad \quad\textit{if} \quad d \quad\textit{is even}\quad .
\end{array} \right.
\end{displaymath}

\noindent \emph{Then, for $n$ satisfying $\gamma^{(2)} =  (2d\,$-$1)(2n$-$2)+3$, and assuming $\gamma:= (2d\,$-$1)\mu + 2\varepsilon$ belongs to $\mathbb{N}^4$ (as well as $\gamma^{(1)}\leq \textit{\textbf{p}}(2d\,$-$\,1)$, if $\textit{\textbf{p}}\geq3$), the divisor $\Lambda(n,d,1,\gamma)$ is nef.}\\

\textbf{Proof.} One only needs to check (straightforward verification!), that any such $\varepsilon$ satisfies \textbf{5.12.}(1),(2) \& (3). $\blacksquare$\\
 
\textbf{Proposition 5.14.}\\
\emph{Let us fix $d \geq 3$  and $\mu\in \mathbb{N}^4$ such that $\mu_o+1\equiv
\mu_j \textnormal{(mod.2)}$ (for $j=1,2,3$), and let $k$ denote the residue \textnormal{(mod.4)} of $\,d+1$. Choose any integer vector $\varepsilon \in \mathbb{Z}^4$ subject to the conditions  $$4\varepsilon^{(2)}= 3+(2d\,\textnormal{-}\,1)(d\,\textnormal{-}\,2+k)\quad\textnormal{and}\quad\gamma :=(2d\,\textnormal{-}\,1)\mu + 2\varepsilon \in \mathbb{N}^4\, ,$$ and let  $n$ satisfy $\gamma^{(2)} =  (2d\,$-$\,1)(2n\,$-$\,2)+3$. Then  $\Lambda(n,d,1,\gamma)$ is not nef}.\\

\textbf{Proof.} Take any vector $\varepsilon \in \mathbb{Z}^4$  satisfying $ \varepsilon^{(2)}=8h^2+3(2k\,$-$\,3)h+k^2\,$-$\,3k +3$.  A straightforward verification shows that $\varepsilon_i^2 \leq \varepsilon^{(2)}<(2d\,$-$\,1)^2, \forall i =0,..,3$ and $4\varepsilon^{(2)}= 3+(2d\,$-$\,1)(d\,$-$\,2+k)$. In particular, $4\varepsilon^{(2)}< 3+(2d\,$-$\,1)^2=4d^2\,$-$\,4d+4$, hence $\varepsilon$ does not satisfy property \textbf{5.12.}(1). Therefore, choosing any $\mu\in \mathbb{N}^4$ such that $\mu_o+1\equiv
\mu_j \textnormal{(mod.2)}$ (for $j=1,2,3$), and defining $\gamma \in \mathbb{N}^4 $ and $n \in\mathbb{N}$ by $\gamma:= (2d\,$-$1)\mu + 2\varepsilon$ and $\gamma^{(2)} =  (2d\,$-$1)(2n$-$2)+3$, respectively, the corresponding divisor $\Lambda(n,d,1,\gamma)$ is not \textit{nef}. $\blacksquare$\\

\textbf{Lemma 5.15.}\\
\emph{Let $(n,d,\gamma) \in \mathbb{N}^*\times \mathbb{N}^*\times \mathbb{N}^4$ be such that $d\geq2$,  $\gamma^{(2)}=(2d\,$-$\,1)(2n\,$-$\,2)+3$ and $\Lambda(n,d,1,\gamma)$ is nef. Then, for any $j=1,2,3$, there exists at most one exceptional curve of the first kind $\widetilde \Gamma \subset \widetilde S$, such that $\widetilde \Gamma\cdot \Lambda(n,d,1,\gamma)=0$ and $\widetilde \Gamma\cdot \widetilde s_j=1$. In particular, the sum of the latter exceptional curves, denoted by  $\widetilde{Z}(n,d,1,\gamma) $, is a reduced divisor with (at most) three irreducible components}. \\

\textbf{Proof.} Straightforward verification again!. $\blacksquare$\\

\textbf{Remark 5.16.}\\
According to Brian Harbourne's results on anticanonical rational surfaces (cf. \cite{Har}), for any \emph{nef} divisor $D \in Pic(\widetilde S)$, such that $\,\textnormal{-}\,\widetilde K \cdot D\geq 2$, the complete linear system $|D|$ is base point free and $dim|D|= \frac{1}{2}D\cdot (D\,$-$\,\widetilde K)$. The following result is in order.\\

\textbf{Lemma 5.17.}\\
\textit{Let  $(n,d,\gamma) \in \mathbb{N}^*\times \mathbb{N}^*\times \mathbb{N}^4$  be such that $d\geq 2 $, $\gamma^{(2)}=(2d\,$-$\,1)(2n\,$-$\,2)+3$, and let  $\Lambda$ and $\widetilde{Z}$ denote, respectively, $\Lambda(n,d,1,\gamma)$ and $\widetilde{Z}(n,d,1,\gamma)$, the divisors defined in \textnormal{\textbf{5.15.}}. Then, $\,\Lambda$ nef implies}:\\
 
\begin{enumerate}
\item $\Lambda\,$-$\,\widetilde C_o\,$-$\, \sum_{j=1}^3 \widetilde s_j \,\,$-$\,\widetilde{Z}$ \textit{is nef};\\
\item $\big|\Lambda \,$-$\,\widetilde C_o \,$-$\,\sum^{3}_{j=1}\widetilde s_j\,$-$\,\widetilde{Z}\big|$ \textit{is base point free};\\

\item $\big|\Lambda\,$-$\,\widetilde C_o\big|=\sum_{j=1}^{3} \widetilde s_j +\widetilde{Z}+ \big|\Lambda \,$-$\,\widetilde C_o \,$-$\,\sum^{3}_{j=1}\widetilde s_j\,\,$-$\,\widetilde{Z}\big|\,$;\\
 
\item $\,\,dim\big|\Lambda\big|=2d\,$-$\,2 $,  $\quad\,dim \big|\Lambda\,$-$\, \widetilde C_o\big|= d\,$-$\,2 \quad $ \textit{and} $\quad h^1\big(\widetilde S,O_{\widetilde S}(\Lambda \,$-$\,\widetilde C_o )\big)=0\,$.\\

\end{enumerate}

\textbf{Definition 5.18.}\\
\emph{Let  $\widetilde p_o \in \widetilde S$ denote the unique point of intersection $\{\widetilde p_o\}:=\widetilde C_o\cap \widetilde s_o$ and consider any divisor $\Lambda: =\Lambda(n,d,1,\gamma)$ as in }\textbf{5.15.}. 
\emph{We define the following subsets of} $\big|\Lambda\big|$:

\begin{align}
 \big|\Lambda\big|_{\widetilde C_o,\widetilde p_o}: &= \Big\{D\in\big|\Lambda\big|,\quad D \cap \widetilde C_o= \{\widetilde p_o\} \quad or \quad \widetilde C_o \subset D\Big\}\, ;\\
\big|\Lambda\big|_{\widetilde C_o,\widetilde p_o}^{\widetilde s_o} :& =\big|\Lambda\big|_{\widetilde C_o,\widetilde p_o} \bigcap \big( \widetilde s_o+\big|\Lambda\,\textnormal{-}\,\widetilde s_o\big|\,\big)   \,.
\end{align}

\textbf{Proposition 5.19.}\\
\textit{If $\Lambda: =\Lambda(n,d,1,\gamma) $ is nef, then}:\\
\begin{enumerate}
\item $\big|\Lambda\big|_{\widetilde C_o,\widetilde p_o}$ \textit{is a $(d\,$-$\,1)$-dimensional subspace of }$\big|\Lambda\big|$;\\

\item $\widetilde C_o $+$\big|\Lambda\,$-$\,\widetilde C_o\big|$ \textit{and $\big|\Lambda\big|_{\widetilde C_o,\widetilde p_o}^{\widetilde s_o}$ are two different hyperplanes of} $\big|\Lambda\big|_{\widetilde C_o,\widetilde p_o}$;\\

\item \textit{any element $\widetilde \Gamma\in \big|\Lambda\big|_{\widetilde C_o,\widetilde p_o}$, in the complement of the latter hyperplanes, is a smooth integral divisor isomorphic to }$\mathbb{P}^1$.
\end{enumerate}
\textbf{Proof.}

(1)  According to \textbf{5.17.}(4), we have $h^1\big(\widetilde S,O_{\widetilde S}(\Lambda \,$-$\,\widetilde C_o )\big)=0$. Hence, the exact sequence of $O_{\widetilde S}$-modules:

\begin{displaymath}
 0\to  O_{\widetilde S} (\Lambda \,\textnormal{-}\,C_o)\to  O_{\widetilde S} (\Lambda )\to O_{\widetilde C_o}(\Lambda ) \to   0\,,
\end{displaymath}  

gives rise to the exact sequence

\begin{displaymath}
 0\to H^0\big( \widetilde S,O_{\widetilde S} (\Lambda \,\textnormal{-}\, \widetilde C_o)\big)\to H^0\big(\widetilde S, O_{\widetilde S}(\Lambda)\big)\to H^0\big(\widetilde C_o, O_{\widetilde C_o}(\Lambda ) \big) \to  0  \,.
\end{displaymath} 

Since $deg\big(O_{\widetilde C_o}(\Lambda ) \big)= d\,$-$\,1$, we can pick a section $f\in H^0\big(\widetilde C_o, O_{\widetilde C_o}(\Lambda ) \big) $ which only vanishes at $\widetilde p_o$ $\big(\,$i.e.: with zero divisor $(f)_o= (d\,$-$\,1)\widetilde p_o \,\big)$, as well as a preimage of $f$, say $v \in H^0\big(\widetilde S, O_{\widetilde S}(\Lambda)\big)$, such that its zero divisor $\widetilde D:= (v)_o \in \big |\Lambda\big |$ only intersects $\widetilde C_o$ at $\widetilde p_o$ $\big(\,$i.e.: $\widetilde D \cap \widetilde C_o= \{ \widetilde p_o\}\,\big)$. Any other section of $O_{\widetilde S}(\Lambda)$, satisfying the same property as $v$, is obtained by adding the image of an arbitrary element of $ H^0\big( \widetilde S,O_{\widetilde S} (\Lambda \textnormal{-}\, \widetilde C_o)\big)$. In other words $\big|\Lambda\big|_{\widetilde C_o,\widetilde p_o}\subset \big|\Lambda\big|$ is the $(d\,$-$\,1)$-dimensional subspace generated by $\widetilde D$ and $\widetilde C_o$+$\big|\Lambda \textnormal{-}\, \widetilde C_o\big|$. 

(2) On the other hand, according to \textbf{5.17.}(2)\&(3), there exists $\widetilde D' \in \big|\Lambda \textnormal{-}\, \widetilde C_o\big|$ avoiding $\widetilde p_o$, in which case $\widetilde C_o$+$\widetilde D' \in \big|\Lambda\big|$ is smooth at $\widetilde p_o$. Up to replacing the former divisor $\widetilde D\in \big|\Lambda\big|$, by the generic element of the pencil generated by $\widetilde D$ and $(\widetilde C_o$+$\widetilde D')$, we can assume hereafter $\widetilde D$ smooth and tangent to $\widetilde C_o$ at $\widetilde p_o$. In particular, for any $\widetilde D''\in \big|\Lambda \textnormal{-}\, \widetilde C_o\big|$, either $\widetilde p_o \notin \widetilde D''$ and $\widetilde C_o$+$\widetilde D''$ is also smooth and tangent to $C_o$ at $\widetilde p_o$, or $\widetilde p_o \in \widetilde D''$ and $\widetilde C_o$+$\widetilde D''$  is singular at $\widetilde p_o$. In both cases, all but one element of the pencil generated by $\widetilde D$ and $\widetilde C_o$+$\widetilde D''$ is smooth and tangent to $C_o$ at $\widetilde p_o$. Therefore, such a generic element is transverse at $\widetilde p_o$ to $\widetilde s_o$, and can not contain $\widetilde s_o$ as an irreducible component. At last, since $\Lambda \cdot \widetilde s_o = 1$, the unique singular element of the latter pencils must belong to $\widetilde s_o$+$\big|\Lambda\,$-$\,\widetilde s_o\big|$. Hence, $\big|\Lambda\big|_{\widetilde C_o,\widetilde p_o}^{\widetilde s_o}$ and $\widetilde C_o $+$\big|\Lambda\,$-$\,\widetilde C_o\big|$ are indeed distinct hyperplanes of $\big|\Lambda\big|_{\widetilde C_o,\widetilde p_o}$. 

(3) Any $\widetilde\Gamma\in \big|\Lambda\big|_{\widetilde C_o,\widetilde p_o}$, in the complement of the latter hyperplanes, has arithmetic genus $0$. Let us also prove its irreducibility. We start remarking that $\widetilde \Gamma$ can only intersect $\widetilde C_o$ at $\widetilde p_o$, and does not contain $\widetilde C_o$ nor $\widetilde s_i \, , (i=0,1,2,3)$, as an irreducible component. Hence, its inverse image $\Gamma^\bot :=\varphi^*(\widetilde \Gamma) \subset S^\bot$ is linearly equivalent to $e^*(nC_o$+$S_o)\,$-$\,s_o^\bot\,$-$\,\sum_{i=o}^3 \gamma_i r_i^\bot$, and neither $C_o^\bot$, nor $s_i^\bot \,(\forall i=0,\dots,3)$, is an irreducible component of $\Gamma^\bot$. In order to check that $\Gamma^\bot$ (hence $\widetilde \Gamma$) is an irreducible curve, by means of \textbf{5.5.}, we still need to show that $r_i^\bot \nsubseteq \Gamma^\bot,  \,\forall i=0,\dots,3$. Otherwise  $\Gamma^\bot$ would have an irreducible component  ${\overline\Gamma}^\bot\subset S^\bot$, linearly equivalent to $e^*(nC_o$+$S_o)\,$-$\,s_o^\bot \,$-$\,\sum_{i=o}^3 {\overline\gamma}_i r_i^\bot$, for some type $\overline\gamma$ strictly bigger than $\gamma$, implying that $\varphi(\,{\overline\Gamma}^\bot) \subset \widetilde S$ has a negative arithmetic genus. Contradiction.! In case $\textit{\textbf{p}} \geq 3$, an analogous line of reasoning shows that $\Gamma^\bot$ can not contain $C_{\textit{\textbf{p}}}^\bot$ as an irreducible component and \textbf{5.5.} still applies.$\blacksquare$\\

Recalling that $MH_X(n,d,1,\gamma)$ is birationally isomorphic to $|\,\Lambda(n,d,1, \gamma)| _{\widetilde C_o,\widetilde p_o}$ (\textbf{5.3.}), we deduce the:\\

\textbf{Corollary 5.20.}\\
\textit{ For any $(n,\mu)\in \mathbb{N}^*\times\mathbb{N}^4$ satisfying $\mu_o+1\equiv \mu_1\equiv \mu_2\equiv \mu_3\textnormal{(mod.2)}$ and $\mu^{(2)} = 2n+1$, (and $\mu^{(1)}\leq \textit{\textbf{p}}$, if $\textit{\textbf{p}}\geq3$), we let $\pi_\mu$ denote the minimal-hyperelliptic $1$-osculating cover associated to the exceptional curve ${\widetilde \Gamma}_\mu \subset \widetilde S$ \textnormal{(cf. \textbf{5.6.} \& \cite{T-V}\S6.2.)}. Then, $\big|\Lambda(n,d,1,\gamma)\big|=\{\widetilde \Gamma_\mu\}$ and $MH_X(n,1,1,\mu)$ reduces to $\{\pi_\mu\}$ }.

\textit{More generally, for any $(n,d,\gamma) \in \mathbb{N}^*\times \mathbb{N}^*\times \mathbb{N}^4$ such that:}
\begin{enumerate}
\item $\gamma_o +1 \equiv \gamma_1\equiv \gamma_2 \equiv \gamma_3(\textnormal{mod}.2) \, \quad$ (and $\gamma^{(1)}\leq \textit{\textbf{p}}$, if $\textit{\textbf{p}}\geq3$), \\
\item $d\geq2 \quad and   \quad \gamma^{(2 )} = (2d\,$-$\,1)(2n\,$-$\,2)+3 \, $, \\
\item $\Lambda(n,d,1,\gamma) \quad is \quad nef,\,$ 
 \end{enumerate}
 
\textit{the moduli space $MH_X(n,d,1,\gamma)$ is birational to $\big |\Lambda(n,d,1,\gamma)\big|_{\widetilde C_o,\widetilde p_o}$}. 

\textit{In particular, $dim\big(MH_X(n,d,1,\gamma)\big)=d\,$-$\,1$, for any $(n,d,\gamma)$ as in} \textbf{5.13.}.\\

 At last, we propose a less conceptual but more geometrical construction of $MH_X(n,d,1,\gamma)$. We will construct $d$ effective divisors $\big\{G^\bot, F_j^\bot, j=0,..,d\,$-$\,2 \big \}$ of $S^\bot$, with birational models given by explicit equations in $\mathbb{P}^1\times X$, which generate all $MH_X(n,d,1,\gamma)$. Hence, any element of $MH_X(n,d,1,\gamma)$ is birational to the zero set of a linear combination of $d$ specific degree-$n$ polynomials with coefficients in \textit{K(X)}, the field of meromorphic functions on $X$.\\

 \textbf{Theorem 5.21.}\\
\noindent   \textit{For any $(n,d,\gamma) \in \mathbb{N}^*\times \mathbb{N}^*\times \mathbb{N}^4$ as in \textnormal{\textbf{5.13.}}, $\big|\,e^*\big(nC_o$}+\textit{$\,(2d$-$1)S_o\big)$-$\,s_o^\bot$-$ \sum _i \gamma_i r_i^\bot  \big|$ contains a
$(d$-$1)$-dimensional subspace with a generic element, say $\Gamma^\perp$, satisfying}:

\begin{enumerate}
\item  $\Gamma^\perp$ \textit{is a} $\tau^\perp$\textit{-invariant smooth irreducible curve of genus }$g
:\ =\frac{1}{2}($-$1$+$\gamma ^{(1)})$;\\

\item \textit{$\Gamma^\bot$ can only intersect $C_o^\bot$ at }$p_o^\bot:=C_o^\bot \cap s_o^\bot$;\\

\item $\varphi(\Gamma^\perp) \subset \widetilde S$ \textit{is isomorphic to } $\mathbb P^1$.\\
  
\end{enumerate}

\textbf{Corollary 5.22.}\\
\noindent \textit{Given $(n,d,\gamma) \in \mathbb{N}^*\times \mathbb{N}^*\times \mathbb{N}^4$ as above, the moduli space $MH_X(n,d,1,\gamma)$  \textnormal{(\textbf{5.2.})} has dimension $d\,$-$\,1$, and its generic element is smooth of genus $g
:\ =\frac{1}{2}($-$1$}$+$ \textit{$\gamma ^{(1)})$.}\\

  \textbf{Proof} of \textbf{Theorem 5.21.}.

  We will only  work out the case $\gamma:= (2d\,$-$1)\mu + 2\varepsilon$, with $\,\,\varepsilon = (0,\,d\,$-$\,1,\,
d\,$-$\,1, \,d\,$-$\,1)\,\,$. 
\indent For any other choice of $\varepsilon$, the corresponding proof runs along the same lines and will be skipped. In our case, the arithmetic genus $g$ and the degree $n$ satisfy: 
\begin{displaymath}
2g+1= (2d\,\textrm{-}\,1)\mu^{(1)}+6(d\,\textrm{-}\,1)\quad \textrm{and} \quad 2n=(2d \,\textrm{-}1) \mu ^{(2)}
\textrm{+} \,4(d\,\textrm{-}1)(\mu_1 \textrm{+} \mu_2 \textrm{+} \mu_3)\textrm{+}\,6d\,\textrm{-}7.
\end{displaymath}

\indent Consider $\overline \mu:\ =\mu$+$\,(1,1,1,1),\,\, \mu' :\  =
\mu$+$\,(0,2,1,1),\,\, \mu''=  \mu$+$\,(0,0,1,1)$, and let $\overline{Z}^{\,\bot}
,{Z'}^\bot,{Z''} ^\bot \subset S ^\bot $ denote the unique $\tau ^\bot$-invariant curves
linearly equivalent to:\\

\noindent 1) $\overline{Z}^{\,\bot} \,\,\,\sim \,e^*(\overline{m}\,C_o$+$\,S_o)\,$-$\,s_o^\bot\,$-$\sum
_i\overline{\mu}_i r_i ^\bot $, \,\,where $\,\,2\overline{m} \,\,$+$ 1 = \overline {\mu}^{\,(2)}$;

\noindent 2) ${Z'}^\bot \,\,\sim \,e^*(m'C_o$+$\,S_1)\,$-$\,s_1^\bot$-$\sum _i\mu_i' 
r_i^\bot$ , \,where \,\,$2m'\,$+$1={\mu'}^{(2)}$;

\noindent 3) ${Z''}^\bot \sim \,e^*(m''C_o$+$\,S_1)\,$-$\,s_1^\bot$-$\sum
_i\mu_i'' r_i^\bot$ , where $\,2m''$+$1={\mu''}^{(2)}$.\\

\indent Moreover, if $\mu_o\neq 0$ we choose $\underline {\mu}\  = \mu $+$\,
($-$\,1,1,1,1)$ and $2\underline {m}\,$+$1= \underline {\mu}^{(2)}$, and let $\underline {Z}^\bot\subset S ^\bot $ denote the unique  $\tau ^\bot$-invariant curve $ \underline {Z}^\bot \sim e^*(\underline {m}C_o$+$\,S_o)\,$-$\,s_o^\bot \,$-$\sum _i\underline {\mu}_{i} r_i^\bot$.

\indent However, if $\mu_o=0$ we will simply put $\underline {Z}^\bot :\  = \overline{Z}^{\,\bot}
$+$ \,2r_o^\bot$, so that in both cases,  
the divisors $D^\bot_0 : = \overline{Z}^{\,\bot}  $+$  \underline {Z}^\bot $+$ \,2s^\bot _0$ and $
D^\bot_1 : = {Z'}^\bot $+$ {Z''}^\bot $+$ \,2s_1^\bot$  will be linearly equivalent. 
Let us also define,  
\begin{align}
\mu _{(1)} :\  = \mu''=\mu + (0,0,1,1), \notag \\
\mu_{(2)}:\ = \mu + (0,1,0,1),\notag \\
\mu_{(3)}: \ = \mu + (0,1,1,0),\notag
\end{align}
  and let $ Z^\bot_{(k)}  (k=1,2,3)$ be the $\tau^\bot$-invariant curve of
$S^\bot $, linearly equivalent to $e^*(m_{(k)}C_o$+$\,S_k)\,$-$\,s_k^\bot $-$
\sum_i\mu_{(k)i}r_i^\bot$, where $ 2m_{(k)}$+$ 1= \sum_i\mu_{(k)i}^2$.

\indent At last, consider $Z^\bot \sim e^*(mC_o$+$S_o) \,$-$\,
s^\bot_o\,$-$\, \sum_i\mu_ir_i^\bot$, where $2m$+$1 = \sum_i \mu_i^2$ (\textbf{5.2.}). Let $\Lambda \in Pic(\widetilde S)$ denote the unique  class such that $\big|\,e^*\big(nC_o$+$\,(2d$-$1)S_o\big)$-$\,s_o^\bot$-$ \sum _i \gamma_i r_i^\bot  \big|=\big|\varphi^*(\Lambda)\big|$. 
The $(d$-$1)$-dimensional subspace of $\big|\varphi^*(\Lambda)\big|$ we are looking for, will be made of all above curves.
 We first remark the following facts :\\

\noindent a) we can check via the adjunction formula, that the divisors $\varphi^*(\Lambda)$  and $\Lambda$ have arithmetic genus $g:\ =\frac{1}{2}($-$1$+$\,\gamma^{(1)})$ and $0$, respectively, and that $ \varphi^*\big(\big|\Lambda\big|\big)$ is equal to $\big|\varphi^*(\Lambda)\big|^{\tau^\bot}$, the sub-space of $\tau^\bot$-invariant 
elements of $\big|\varphi^*(\Lambda)\big|$;\\

\noindent b) the  $d\,$-$\,1$ divisors 

\begin{displaymath}
F^\bot_j : = C^\bot_o \textrm{+} \sum _{k=1}^{3}(Z^\bot_{(k)}\textrm{+}2s^\bot_k) \textrm{+}
jD^\bot_o \textrm{+} \,(d\,\textrm{-}2\,\textrm{-}\,j)D^\bot_1, \,\,\quad \quad j=0,...,d\,\textrm{-}2,
 \end{displaymath}
as well as 
\begin{displaymath}  
 G^ \bot : = Z^\bot \textrm{+} \,(d\,\textrm{-}1)D^\bot_o, 
 \end{displaymath}
are $\tau^\bot$-invariant, belong to $\big|\varphi^*(\Lambda)\big|$ and have $p^\bot_o :  = C^\bot_o
\cap s_o^\bot$ as their unique common point;\\

\noindent c) the curve $F^\bot_o$ is smooth at $p^\bot_o$, while any other $F^\bot_j$ has multiplicity $1<2j+1< 2d$ at $p^\bot_o$. In particular, they span a $(d\,$-$\,2)$-dimensional subspace of $\big|\varphi^*(\Lambda)\big|$, having a generic element smooth and  
transverse to $s_o^\bot$ at $p^\bot_o$;\\

\noindent d) the curve $G^\bot$ has multiplicity $2d$ at $p_o^\bot$, and no common irreducible component with any $F^\bot_j (\,\forall j=0, \ldots,d\,$-$\,2)$, implying that $\langle \,G^\bot,\,F_j^\bot,\,j=0,..,d\,$-$\,2\,\rangle $, the $(d\,$-$\,1)$-dimensional subspace they span in $\big|\varphi^*(\Lambda)\big|$, is component-free;\\

\noindent e) any irreducible curve $\Gamma^\bot \in \langle \,G^\bot,\,F_j^\bot,\,j=0,..,d\,$-$\,2\,\rangle$ projects onto a smooth irreducible curve (isomorphic to $\mathbb{P}^1$). In particular $\Gamma^\bot$ must be smooth outside $\cup_{i=0}^3r_i^\bot$.\\

\noindent f) the curves $G^\bot$ and $F_o^\bot$ have no common point on any $r_i^\bot$ ($i=0,..,3$), implying that $\Gamma^\bot$, the generic element of $\langle \,G^\bot,\,F_j^\bot,\,j=0,..,d\,$-$\,2\,\rangle$, is smooth at any point of $\cup_{i=0}^3r_i^\bot$ and satisfies the announced properties, i.e.:\\

(1) $\Gamma^\bot$ is $\tau^\bot$-invariant, smooth and satisfies the irreducibility criterion \textbf{5.5.};\\

(2) $p_o^\bot$ is the unique base point of the linear system and $\Gamma^\bot \cap C_o^\bot= \{p_o^\bot\}$;\\

(3) its image $\varphi(\Gamma^\bot) \subset \widetilde S$ is irreducible, linearly equivalent to $ \Lambda(n,d,1,\gamma)$ and 
\indent of arithmetic genus
$\frac{1}{4}\big((2d\,$-$1)(2n$-$2)$+$\,3\,$-$\, \gamma^{(2)}\big)  = 0$; hence, isomorphic to $\mathbb P^1$. $\blacksquare $\\

\textbf{Proof} of \textbf{Corollary 5.22.}.

 The degree-$2$ projection $\varphi: \Gamma^\bot \longrightarrow \varphi(\Gamma ^\bot)$ is ramified at $p_o^\bot$ and $\varphi(\Gamma^\bot)$ is isomorphic to $\mathbb{P}^1$. Moreover, $\Gamma^\bot$ is a smooth irreducible curve linearly equivalent to $ \varphi^*\big(\Lambda(n,d,1,\gamma)\big)$, of  arithmetic genus 
$g:=\frac{1}{2}(\gamma^{(1)}\,$-$\,1) $. 

In other words, the natural projection $ (\Gamma ^\bot,\,p_o^\bot) \subset (S^\bot ,\,p_o^\bot)\stackrel{\pi_{S^\bot}}{\longrightarrow} (X,q)$ is a smooth degree-$n$ \textit{minimal-hyperelliptic $d$-osculating cover} of type $\gamma $,
and genus $g$, such that $(2n$-$2)(2d\,$-1$)$+$\,3=\gamma ^ {(2)} $ and
$2g\,$+$1=\gamma^{(1)}.$ $  \blacksquare $\\

\end{document}